\documentclass[10pt]{article}

\usepackage{algorithmic}
\usepackage{algorithm}
\usepackage{enumerate}
\usepackage{fullpage,graphicx,amssymb,amsmath,amsthm}
\usepackage{hyperref}
\usepackage{subfigure}
\usepackage[usenames,dvipsnames]{color}
\usepackage{multirow}
\usepackage{booktabs}

\newcommand{\be}{\begin{enumerate}}
\newcommand{\ee}{\end{enumerate}}

\newtheorem{corollary}{Corollary}

\newtheorem{proposition}{Proposition}

\newcommand{\cl}{{\rm cl}} % Closure
\newcommand{\clconv}{{\rm clconv}} % Closed convex hull
\newcommand{\D}{\mathbb{D}} % Distance metric
\renewcommand{\Re}{\mathbb{R}} % Real numbers
\newcommand{\sign}{\ensuremath{\text{sign}}}
\newcommand{\Hrel}{\text{Hrel}} % Hull relaxation
\newcommand{\LD}{\text{LD}} % Lagrangian dual

\newcommand{\vlambda}{\boldsymbol{\lambda}}
\newcommand{\vmu}{\boldsymbol{\mu}}
\newcommand{\vnu}{\boldsymbol{\nu}}

\renewcommand{\v}[1]{\ensuremath{\mathbf{#1}}}
\newcommand{\mc}{\mathcal}

\def\st{{\rm s.t.}}

\theoremstyle{remark}
\newtheorem{remark}{Remark}

\begin{document}

\title{Pseudo basic steps:\\ Bound improvement guarantees from Lagrangian decomposition\\ in convex disjunctive programming}
\author{Dimitri J.\ Papageorgiou and Francisco Trespalacios  \\
{\small Corporate Strategic Research}\\
{\small ExxonMobil Research and Engineering Company}\\
{\small 1545 Route 22 East, Annandale, NJ 08801 USA}\\
{\small \texttt{\{dimitri.j.papageorgiou,francisco.trespalacios\}@exxonmobil.com}} \\
}

\maketitle

\begin{abstract}
An elementary, but fundamental, operation in disjunctive programming is a basic step,
which is the intersection of two disjunctions to form a new disjunction. Basic steps bring a disjunctive set in regular form closer to its disjunctive normal form and, in turn, produce {\color{black}relaxations that are at least as tight}.  An open question is: What are guaranteed bounds on the improvement from a basic step? In this paper, using properties of a convex disjunctive program's hull reformulation and multipliers from Lagrangian decomposition, we introduce an operation called a \textit{pseudo} basic step and use it to provide provable bounds on this improvement along with techniques to exploit this information when solving a disjunctive program as a convex MINLP.
Numerical examples illustrate the practical benefits of these bounds. 
In particular, on a set of $K$-means clustering instances, we make significant bound improvements relative to state-of-the-art commercial mixed-integer programming solvers.

\vspace{2 mm}
\end{abstract}

\section{Introduction}

Disjunctive programming is optimization over disjunctive sets, 
i.e., sets of inequalities connected to one another by the operations of conjunction or disjunction.
Conceived by Balas in the early 1970s, disjunctive programming was originally advertised as ``a class of problems which subsumes pure and mixed-integer programs and many other nonconvex programming problems'' \cite{balas1998disjunctive}.  
Within the integer programming community, Balas's seminal ideas have gone on to influence a number of topics ranging from lift-and-project cuts \cite{balas:lift_and_project:1993} and disjunctive cuts to intersection cuts for convex lattice-free sets and more.
For an enjoyable first-hand account of disjunctive programming, the reader is encouraged to peruse \cite[Chapter 10]{balas:2010}.

In the process systems engineering community, which has also made numerous contributions to the field of optimization,
%another branch on the disjunctive programming tree was formed 
another branch of disjunctive programming was formed that goes by the name \textit{generalized disjunctive programming} (GDP).
Raman and Grossmann \cite{raman1994modelling} generalized Balas's initial paradigm in two important ways.
First, whereas Balas focused on a purely linear setting, often referring to disjunctive programming as ``optimization over a union of polyhedra'' or ``linear programming with disjunctive constraints,'' Raman and Grossmann extended this framework to include arbitrary nonconvex sets.
Second, and perhaps more importantly, they introduced logic propositions in terms of Boolean variables to provide 
%a more expressive modeling paradigm than mixed-integer programming.
``an alternative modeling framework to mixed-integer linear/nonlinear programming that is more adept at translating physical intuition of engineering principles into rigorous mathematical formalism'' \cite{sawaya:hierarchy:2012}. 
To be precise, a student of mixed-integer programming (MIP) is taught to write down a model expressed solely in terms of two primitives: continuous/integer decision variables, which are either real- or integer-valued, and algebraic constraints.  Meanwhile, a student of GDP is taught to model with, in addition to the aforementioned primitives, disjunctions (including implications) and Boolean variables.  With these additions, GDP is arguably a more expressive modeling paradigm than mixed-integer programming. For a more detailed discussion on GDP modeling, see Grossmann and Trespalacios \cite{grossmann2013systematic}.  
%For example, whereas a student of mixed-integer programming is taught to write down a model expressed solely in terms of continuous and integer decision variables and algebraic constraints, a student of GDP is taught to 
In this paper, we do not consider Boolean variables and logic propositions as found in GDP; hence, we omit the term ``generalized.''

%Therefore, it is relevant to mention some past work on convex disjunctive programming.
There are several other noteworthy papers on or involving convex disjunctive programming.
Ceria and Soares \cite{ceria1999convex} showed how to construct the convex hull of the union of a finite number of convex sets by exploiting the perspective function.  Stubbs and Mehrotra \cite{stubbs1999branch} devised a branch-and-cut algorithm for 0-1 mixed-integer convex optimization problems.
Ruiz and Grossmann \cite{ruiz:hierarchy:2012} extended hierarchies of disjunctive programming and linear GDP in \cite{balas:hierarchy:1985} and \cite{sawaya:hierarchy:2012}, respectively, to nonlinear convex GDPs.
Trespalacios and Grossmann \cite{trespalacios2014algorithmic} proposed an algorithmic approach to improve convex GDP formulations that involves, among other steps, the iterative application of basic steps and the judicious selection of where to apply them.
Our work can be viewed as an extension of the Lagrangian relaxation introduced by Trespalacios and Grossmann \cite{trespalacios:lagrangean:2016} for \textit{linear} GDPs.
Bonami et al. \cite{bonami:2015} and Belotti et al. \cite{belotti:2016} present techniques for handling indicator constraints, i.e., constraints that either hold or are relaxed depending on the value of a binary variable.  Since an implication is logically equivalent to a disjunction, many of the techniques suggested make use of disjunctive programming concepts. 
The findings in this paper are meant, in part, to support the claim made in 
Belotti et al. \cite{belotti:2016} ``... that aggressive bound tightening is often overlooked in MIP, while it represents a significant building block for enhancing MIP technology when indicator constraints and disjunctive terms are present.''

%Although a GDP is often expressed with both logical propositions and Boolean variables,
%for our purposes, we will simplify notation and express a nonlinear GDP in the following compact form:
Consider a nonlinear disjunctive program 
\begin{subequations} \label{model:nonlinear_gdp_set_based}
\begin{alignat}{4}
\min~~&  f(\v{x}) & & \\
\st~~& g(\v{x}) \leq \v{0} & & \\
& \v{x} \in \mc{F} = \bigcap_{k \in \mc{K}} \bigcup_{i \in \mc{D}_k} \mc{S}_{ki} & & ~,
\end{alignat}
\end{subequations}
where  $f$ is a nonlinear function, $g(\v{x}) \leq \v{0}$ denotes a system of ``global'' constraints that are independent of the disjunctions, $\mc{K} = \{1,\ldots,K\}$ is a set of disjunctions, 
$\mc{D}_k$ is the set of disjuncts associated with the $k$th disjunction, 
and $\mc{S}_{ki} \subset \Re^n$ is a closed and bounded set describing the feasible region of the $i$th disjunct in the $k$th disjunction.
In this paper, we focus on the special case when $f$ is a convex function and all disjuncts are convex, which can be written as: 
\begin{subequations} \label{model:convex_gdp_set_based}
\begin{alignat}{4}
\min~~& \v{c}^\top \v{x} & & \\
\st~~& \v{x} \in \bigcap_{k \in \mc{K}} \bigcup_{i \in \mc{D}_k} \mc{C}_{ki}~. & & 
\end{alignat}
\end{subequations}
Notice that in \eqref{model:convex_gdp_set_based}, we consider a linear objective function since one can always minimize an auxiliary variable $z$ and include a constraint $z \geq f(\v{x})$. 
Likewise, there are no global constraints $g(\v{x}) \leq \v{0}$.
This is without loss of generality since global constraints can always be handled in at least two alternate ways: they can be (1) treated as a single (improper) disjunction possessing a single disjunct, or (2) embedded within each disjunct,
although this obviously comes at the expense of working with disjuncts having more inequalities. 

To obtain a disjunctive program's so-called \textit{hull reformulation}, it is convenient{\color{black}, but not necessary,} to work with an extended formulation in which ``copies'' (also known as \textit{disaggregated variables}) of the original decision variables are created, one copy for each disjunct.  This leads to the following formulation:  
\begin{subequations} \label{model:disaggregated_convex_gdp_set_based}
\begin{alignat}{4}
z^* = 
\min~~& \v{c}^\top \v{x} & & \\
\st~~& \v{x} = \v{v}_k  & & \qquad \forall~k \in \mc{K} \\
& \v{v}_k \in \mc{F}_k = \bigcup_{i \in \mc{D}_k} \mc{C}_{ki} & & \qquad \forall~k \in \mc{K} \label{eq:F_k}
\end{alignat}
\end{subequations}
The hull relaxation can then be expressed as
\begin{subequations} \label{model:hr_disaggregated_convex_gdp_set_based}
\begin{alignat}{4}
(HR_{\mc{K}}) \quad z_{\mc{K}}^{\Hrel} = 
\min~~& \v{c}^\top \v{x} & & \\
\st~~& \v{x} = \v{v}_k  & & \qquad \forall~k \in \mc{K} \quad (\vlambda_k \in \Re^n) \label{eq:hr_disaggregation_constraints} \\
& \v{v}_k \in \clconv \left( \mc{F}_k \right) & & \qquad \forall~k \in \mc{K} ~,
\end{alignat}
\end{subequations}
where $\clconv$ denotes the convex closure of a set.
Lagrange multipliers, if they exist, that link the original variables $\v{x}$ with their copies $\v{v}_k$ are expressed as $\vlambda_k \in \Re^n$.
%We assume that, when represented algebraically, the Slater condition (see, e.g., Boyd and Vandenberghe \textit{Convex Optimization}) holds to ensure the existence of Lagrange multiplier vectors $\vlambda_k$ for each disjunction $k \in \mc{K}$.  Such a constraint qualification also ensures that strong duality holds.

\begin{remark}
Solving the hull relaxation \eqref{model:hr_disaggregated_convex_gdp_set_based} is a convex optimization problem and, thus, can be done in polynomial time. 
\end{remark}

{\color{black} 
Our results assume that the original convex disjunctive program~\eqref{model:convex_gdp_set_based} has been reformulated so that we can exploit the multipliers $\vlambda_k$ corresponding to the constraints~\eqref{eq:hr_disaggregation_constraints}. 
In addition, the prevailing theory typically assumes an extended formulation using disaggregated variables to obtain a tractable algebraic representation of the hull reformulation \cite{ceria1999convex,ruiz:hierarchy:2012, stubbs1999branch}.

Given the progress in general purpose optimization solvers, disjunctive programs are typically re-formulated as mixed-integer nonlinear programs (MINLPs) when solved numerically. To achieve this, each disjunct is encoded with an algebraic description of the form $\mc{C}_{ki} = \{ \v{x} \in \Re^n : g_{ki}(\v{x}) \leq \v{0} \}$, where $g_{ki}$ are convex functions for all $k \in \mc{K}$ and $i \in \mc{D}_k$.
As shown in Ruiz and Grossmann \cite{ruiz:hierarchy:2012}, based on the work of \cite{ceria1999convex} and \cite{stubbs1999branch}, a convex disjunctive program can be expressed as the following convex MINLP: 
\begin{subequations} \label{model:hr_disaggregated_convex_gdp_algebraic}
\begin{alignat}{4}
\min_{\v{v},\v{x},\v{y},\vnu}~~& \v{c}^\top \v{x} & & \\
\st~~~~& \v{x} = \v{v}_k  & & \qquad \forall~k \in \mc{K} \quad (\vlambda_k \in \Re^n) \\
& \v{v}_k = \sum_{i \in \mc{D}_k} \vnu_{ki}  & & \qquad \forall~k \in \mc{K} \label{eq:v_eq_nu} \\
& (\cl~g_{ki}')(\vnu_{ki},y_{ki}) & & \qquad \forall~k \in \mc{K}, i \in \mc{D}_k \\
& \sum_{i \in \mc{D}_k} y_{ki} = 1 & & \qquad \forall~k \in \mc{K} \\
& -\v{L}y_{ki} \leq \vnu_{ki} \leq \v{L}y_{ki} & & \qquad \forall~k \in \mc{K}, i \in \mc{D}_k \\
& y_{ki} \in \{0,1\} & & \qquad \forall~k \in \mc{K}, i \in \mc{D}_k~, \label{eq:y_bounds}
\end{alignat}
\end{subequations} 
where generically $(\cl~g')(\v{x},y)$ denotes the closure of the perspective function $g'$ of $g(\v{x})$ at $(\v{x},y)$ and is given by $g'(\v{x},y) = yg(\v{x}/y)$.
In particular, $(\cl~g_{ki}')(\vnu_{ki},y_{ki}) = y_{ki}g_{ki}(\vnu_{ki}/y_{ki})$ if $y_{ki} > 0$ and $(\cl~g_{ki}')(\vnu_{ki},y_{ki}) = 0$ if $y_{ki} = 0$.
Here, $\v{L}$ denotes a vector of finite bounds on the absolute value of $\v{x}$, which exists due to the assumption that all disjuncts are compact.

For a general convex disjunctive program, attempting to solve hull formulation~\eqref{model:hr_disaggregated_convex_gdp_algebraic} can lead to numerical difficulties when employing a branch-and-cut algorithm for at least two reasons.  First, the resulting convex MINLP is large leading to large, possibly ill-conditioned, convex subproblems at each search tree node. Second, the functions $(\cl~g_{ki}')(\vnu_{ki},y_{ki})$ are not differentiable at $y_{ik} = 0$.
A popular approach to combat these non-differentiabilities was introduced in Sawaya~\cite{sawaya2006reformulations} and relies on approximating the functions $(\cl~g_{ki}')$ in such a way that the formulation is tight when $y_{ik} \in \{0,1\}$, but slightly relaxed at non-binary solutions. 

When working with specially structured convex disjunctive programs, many of these closure and differentiability concerns disappear. For linear disjunctive programs (i.e., all disjuncts $\mc{C}_{ki}$ are polytopes), Balas \cite{balas:1979} provides an explicit extended formulation for the hull reformulation.
For conic quadratic disjunctive programs, each disjunct is the intersection of one or more conic quadratic constraints $\v{x}^\top \v{Q}\v{x} + \v{b}^\top \v{x} + \gamma \leq 0$.  Omitting subscripts, consider a disjunct with a single conic quadratic constraint.
The set corresponding to the perspective function of a conic quadratic function is given by
$\left\{ (\v{x},y) \in \Re^n \times [0,1] : \v{x}^\top \v{Q}\v{x}/y + \v{b}^\top \v{x} + \gamma y \leq 0 \right\}$.
%$\mc{Q}' = \left\{ (\v{x},y) \in \Re^n \times [0,1] : \frac{\v{x}^\top \v{Q}\v{x}}{y} + \v{b}^\top \v{x} + \gamma \leq 0 \right\}$.
Although this set is not closed at $y=0$ since $\v{x}^\top \v{Q}\v{x}/y$ is not defined at $y=0$,
it is possible to lift this set into a higher dimensional space and eliminate this issue:
$\left\{ (\v{x},y,t) \in \Re^n \times [0,1] \times \Re_+ : t + \v{b}^\top \v{x} + \gamma y \leq 0 , \v{x}^\top \v{Q}\v{x} \leq ty \right\}$.
This set is conic quadratic representable as the constraint $\v{x}^\top\v{Q}\v{x} \leq ty$ is a rotated second-order cone constraint (Ben-Tal and Nemirovski~\cite{bental2001lectures}).
We will use this result multiple times in our computational experiments.
Finally, we have deliberately chosen to work with the consolidated formulation~\eqref{model:disaggregated_convex_gdp_set_based},
in which constraints~\eqref{eq:F_k} subsume constraints~\eqref{eq:v_eq_nu}-\eqref{eq:y_bounds},
as the bulk of our results do not require the latter representation. 

}

\section{Basic Steps} \label{sec:basic_steps}

%Basics steps are not in the MIP repertoire.
Basics steps are not a commonly used tool in a mixed-integer programmer's toolbox.
A quick look in any of the major solvers' API reveals that a basic step is not even a callable operation. 
Thus, it is natural to ask: Why should one care about basic steps in the first place?
We attempt to provide some answers to this question and revive interest in this basic operation.

As described by Balas \cite{balas:hierarchy:1985},
a disjunctive set $\mc{F}$ can be expressed in many different forms that are logically equivalent 
and can be obtained from each other by considering $\mc{F}$ as a logical expression whose statement forms are inequalities,
and applying the rules of propositional calculus.
Among these equivalent forms, the two extremes are the \textit{conjunctive normal form} (CNF)
$$
\mc{F} = \bigcap_{k \in \mc{K}} \mc{E}_k~,
$$
where each $\mc{E}_k$ is an elementary disjunctive set, which, for convex disjunctive programming, is 
the union of finitely many convex inequalities $g_j(\v{x}) \leq 0$, 
%a disjunction in which all disjuncts possess a single convex inequality $g(\v{x}) \leq 0$, 
and the \textit{disjunctive normal form} (DNF)
$$
\mc{F} = \bigcup_{i \in \mc{D}} \mc{C}_{i}~,
$$
where each $\mc{C}_{i}$ is a convex set (and compact, in our setting).
Why is DNF important? Balas \cite[Theorems 3.3 and 3.4]{balas:hierarchy:1985} showed that the convex hull of a linear disjunctive program in DNF is the projection of a higher dimensional polyhedron onto $\Re^n$. 
Practically speaking, this means that solving the hull relaxation of a linear disjunctive program in DNF, which is a linear program (albeit a potentially very large one), yields an optimal solution to the original problem. 
Meanwhile, most discrete optimization problems are stated in the form of an intersection of elementary disjunctions, that is, in CNF.

\begin{figure}[htbp]
\begin{center}
\includegraphics[width=5.0in,height=2.5in]{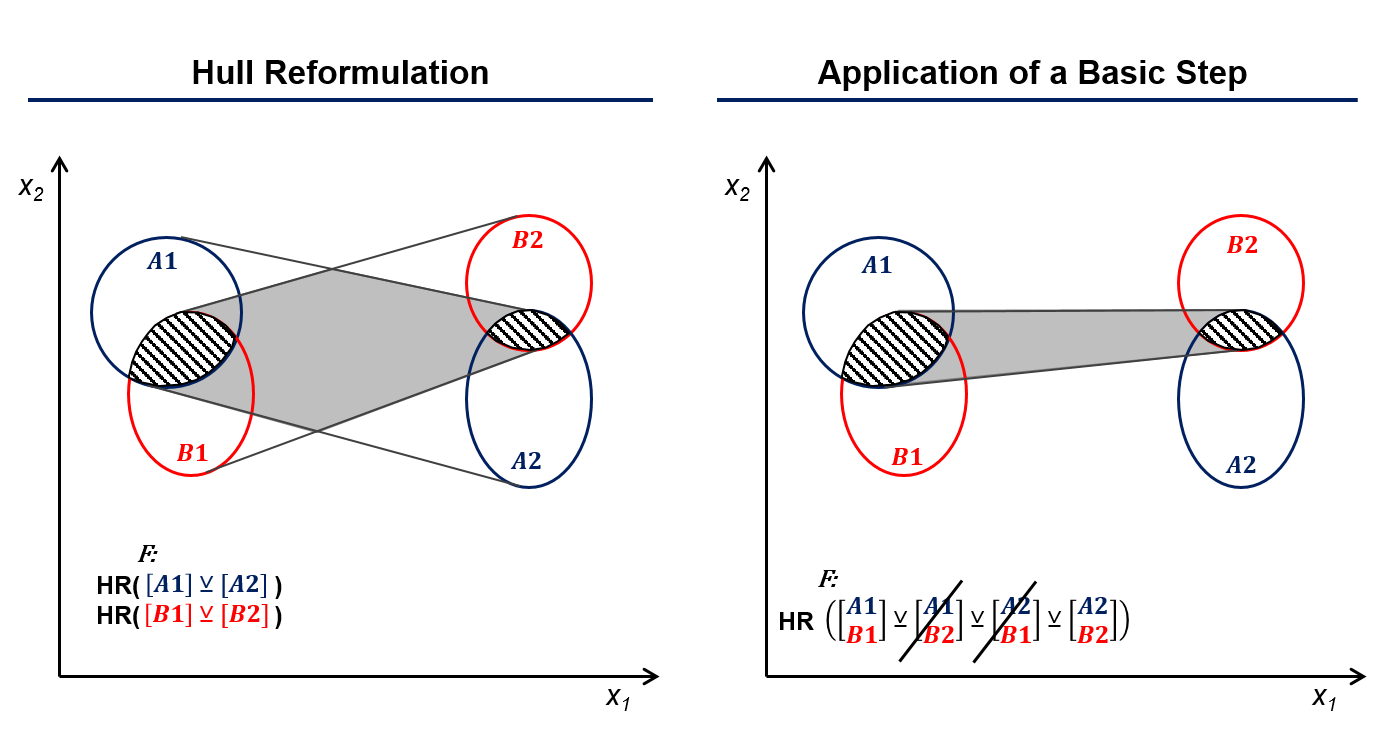}
%\includegraphics{Circles_before_and_after_basic_step.png}
%\vspace{-0.25in}
\caption{Illustration of basic steps.}
\label{fig:basic_step_with_circles}
\end{center}
\end{figure}

A basic step is an elementary, but fundamental, operation in disjunctive programming, tantamount in importance to a branching rule or cutting plane in mixed-integer programming.
In words, a basic step applied to a disjunctive set $\mc{F} = \bigcap_{k \in \mc{K}} \bigcup_{i \in \mc{D}_k} \mc{S}_{ki}$
is the operation of taking two disjunctions $p$ and $q$ in $\mc{K}$ and bringing their intersection into disjunctive normal form.  In set notation, a basic step takes some $p,q \in \mc{K} (p \neq q)$ and brings $\left( \cup_{i \in \mc{D}_p} \mc{S}_{pi} \right) \cap \left( \cup_{j \in \mc{D}_{q}} \mc{S}_{q j} \right) $ into DNF by replacing it with the set
\begin{equation}
\bigcup_{\substack{i \in \mc{D}_p\\j \in \mc{D}_{q}}} (\mc{S}_{pi} \cap \mc{S}_{qj})~.
\end{equation} 

Since every application of a basic step replaces two disjunctions with a single disjunction,
the set $\mc{F}$ can be brought into DNF by $K-1$ applications of basic steps.
%However, one need not bring $\mc{F}$ into DNF in order for basic steps to be useful.
However, just as it is rarely necessary to find the integer hull of a mixed-integer program in order to solve it,
one need not bring a disjunctive program into DNF to obtain useful bounds.  In fact, a few judiciously chosen basic steps may suffice to adequately improve the lower bound.
Figure \ref{fig:basic_step_with_circles} illustrates an application of a basic step between two disjunctions $A$ and $B$, each possessing two disjuncts.  The striped area delineates the true feasible region while the shaded (and striped) area marks the feasible region of the hull relaxation before and after the basic step.  Clearly, the hull relaxation is much tighter after the basic step.
%From the above discussion, it should now be clear why a basic step is not regularly used in mixed-integer programming computations. It typically requires the generation of additional decision variables which reek havoc on the efficiency of most solvers. 
{ \color{black}
Note that much of the gain is due to the fact that, after the basic step, two of the resulting disjuncts are empty sets and thus can be discarded.

From the above discussion, it should be clear that iteratively applying basic steps results in fewer disjunctions, but with certain disjunctions possibly possessing an exponential number of disjuncts. 
If we were to construct hull reformulations for these disjunctive programs as formulated in \eqref{model:hr_disaggregated_convex_gdp_algebraic}, which require disaggregated variables $\vnu_{ki}$ and $y_{ki}$ for each disjunct,
the resulting convex MINLP formulations would likely require the generation of additional decision variables after each basic step.  Adding decision variables during the course of a branch-and-cut algorithm reeks havoc on the efficiency of most solvers, which helps explain why basic steps are not regularly used in mixed-integer programming computations. 
}

Since basic steps have the potential to produce tighter relaxations, 
it is natural to ask: Which basic steps result in the tightest possible relaxation?
Since an answer to this question was elusive at the time,
Balas \cite[Theorem 4.5]{balas:hierarchy:1985} provided an answer to a simpler, but nonetheless important question:
When will a basic step be ineffective?  
His result was specialized to the linear setting, but can be trivially extended to the convex setting below.  
\begin{proposition}
For $k=1,2$, let $\mc{F}_k = \cup_{i \in \mc{D}_k} \mc{C}_{ki}$, where each $\mc{C}_{ki}$ is a compact convex set.  Then 
\begin{equation}
\mc{T}_L := \clconv( \mc{F}_1 \cap \mc{F}_2 ) = \clconv( \mc{F}_1 ) \cap \clconv( \mc{F}_2 ) =: \mc{T}_R 
\end{equation}
if and only if every extreme point of $\mc{T}_R$ is an extreme point of $\mc{C}_{1p} \cap \mc{C}_{2q}$ for some $(p,q) \in \mc{D}_1 \times \mc{D}_2$.
\end{proposition}
%Balas \cite{balas1998disjunctive}, Sawaya and Grossmann ?, Ruiz and Grossmann \cite{ruiz:hierarchy:2012}, and Trespalacios and Grossmann \cite{trespalacios2014algorithmic} provide others ... TODO
Unfortunately, evaluating this property is not practical in general, and can only be done tractably for some specific types of constraints. Meanwhile, for convex disjunctive programs, Ruiz and Grossmann \cite{ruiz:hierarchy:2012} characterized additional types of constraints for which basic steps do not improve relaxations. Trespalacios and Grossmann \cite{trespalacios2014algorithmic} provided various heuristic rules concerning when to apply basic steps based on whether two disjunctions have variables in common, the number of disjuncts in a disjunction, and a parameter they define as the characteristic value of a disjunction, which represents how much the objective function changes if branched on the disjunction. The latter idea is similar to techniques used in strong/pseudocost/reliability branching.

\section{Main results}

In this section, we present our main result by answering the question: What is a bound improvement guarantee from a basic step? 
To accomplish this, we take advantage of a Lagrangian relaxation of a disaggregated convex disjunctive program,
for which the linear case was originally considered in Trespalacios and Grossmann \cite{trespalacios:lagrangean:2016}.  
The review work by Guignard \cite{guignard2003lagrangean} discusses how Lagrangian relaxation can be used in different solution methods and applications.
In Section \ref{sec:main_results}, along with answering the above question, we introduce the central tool of this work, which we call a pseudo basic step.
Section \ref{sec:illustrative_example} then presents an illustrative example as well as insight into why pseudo basic steps may be less productive than actual basic steps.  Discussion to contextualize pseudo basic steps follows in Section \ref{sec:discussion}.

\subsection{Partition relaxation and guaranteed bound improvements from basic steps} \label{sec:main_results}

Consider a partial Lagrangian relaxation of a disaggregated convex disjunctive program \eqref{model:disaggregated_convex_gdp_set_based} described by the set $\mc{K} = \{1,\ldots,K\}$ of disjunctions: 
\begin{subequations} \label{model:lr_disaggregated_convex_gdp_set_based}
\begin{alignat}{4}
LR_{\mc{K}}(\vlambda_1,\ldots,\vlambda_{K}) \quad =~
\min~~& (\v{c} - \sum_{k \in \mc{K}} \vlambda_k )^\top \v{x} + \sum_{k \in \mc{K}} \vlambda_k^\top \v{v}_k & & \\
\st~~& \v{v}_k \in \mc{F}_k = \bigcup_{i \in \mc{D}_k} \mc{C}_{ki} & & \qquad \forall~k \in \mc{K} \label{eq:vk_in_Fk} \\
& \v{x} \in \Re^n~.
\end{alignat}
\end{subequations}
{\color{black}The qualifier ``partial'' serves to emphasize that, while there may be other inequalities in an algebraic representation of the sets $\clconv(\mc{F}_k)$ that could be dualized, we only relax the constraints $\v{x}=\v{v}_k$ linking the original decision vector and its copies.}
Moreover, consider the Lagrangian dual problem 
\begin{equation} \label{model:Lagrangian_dual_problem}
z_{\mc{K}}^{\LD} = \max \left\{ LR_{\mc{K}}(\vlambda_1,\ldots,\vlambda_{K}) : \vlambda_k \in \Re^n~\forall k \in \mc{K} \right\}~.
\end{equation}

\begin{proposition} \label{prop:optimal_multipliers_condition}
If optimal multipliers $\vlambda_1^*,\ldots,\vlambda_{K}^*$ exist for \eqref{model:Lagrangian_dual_problem}, they must satisfy the condition $\sum_{k \in \mc{K}} \vlambda_k^* = \v{c}$.
\end{proposition}
\proof If $\sum_{k \in \mc{K}} \vlambda_k \neq \v{c}$, the minimization problem \eqref{model:lr_disaggregated_convex_gdp_set_based} is unbounded as $\v{x}$ is unrestricted. \qed

Henceforth, we only consider multipliers $\vlambda_1,\ldots,\vlambda_{K}$ such that $\sum_{k \in \mc{K}} \vlambda_k = \v{c}$.  Thus, the Lagrangian relaxation problem \eqref{model:lr_disaggregated_convex_gdp_set_based} can be further decomposed into $K$ subproblems{\color{black}, each one being a convex disjunctive program in DNF,} as follows: 
\begin{equation} \label{model:lr_disaggregated_convex_gdp_set_based_st_optimality_condition}
LR_{\mc{K}}(\vlambda_1,\ldots,\vlambda_{K}) =
\sum_{k \in \mc{K}} \min \left\{ \vlambda_k^\top \v{v} : \v{v} \in \mc{F}_k \right\}~.
\end{equation}
{\color{black}
Re-stating the Lagrangian relaxation as in \eqref{model:lr_disaggregated_convex_gdp_set_based_st_optimality_condition} reveals that both formulations~\eqref{model:lr_disaggregated_convex_gdp_set_based} and \eqref{model:lr_disaggregated_convex_gdp_set_based_st_optimality_condition} are also partial Lagrangian relaxations of the hull relaxation~\eqref{model:hr_disaggregated_convex_gdp_set_based}, a relationship captured in the following proposition. 

\begin{proposition} \label{prop:weak_and_strong_duality} 
The hull relaxation~\eqref{model:hr_disaggregated_convex_gdp_set_based} and the Lagrangian dual problem~\eqref{model:Lagrangian_dual_problem} are weak duals of one another, i.e., $z_{\mc{K}}^{\Hrel} \geq z_{\mc{K}}^{\LD}$. 
Moreover, if optimal multipliers $\vlambda_1^*,\ldots,\vlambda_{K}^*$ exist to \eqref{model:hr_disaggregated_convex_gdp_set_based},
then they are optimal for \eqref{model:Lagrangian_dual_problem} and $z_{\mc{K}}^{\Hrel} = z_{\mc{K}}^{\LD}$.
\end{proposition}
\proof
The weak duality claim follows from the fact that each subproblem $\min \left\{ \vlambda_k^\top \v{v} : \v{v} \in \mc{F}_k \right\}$ in \eqref{model:lr_disaggregated_convex_gdp_set_based_st_optimality_condition} is the minimization of a linear function over the union of compact convex sets, which is equivalent to minimizing the same linear function over the convex closure of the union of these sets (see, e.g., \cite[Section 2.5]{ruiz:hierarchy:2012}). 
Meanwhile, the existence of optimal multipliers $\vlambda_1^*,\ldots,\vlambda_{K}^*$ to the hull relaxation~\eqref{model:hr_disaggregated_convex_gdp_set_based} implies, by \cite[Definition 5.1.1]{bertsekas1999nonlinear}, that $z_{\mc{K}}^{\Hrel} = LR(\vlambda_1^*,\ldots,\vlambda_{K}^*)$. Hence, $z_{\mc{K}}^{\Hrel} = z_{\mc{K}}^{\LD}$ and $\vlambda_1^*,\ldots,\vlambda_{K}^*$ are optimal for \eqref{model:Lagrangian_dual_problem} (see also \cite[Proposition 5.1.4]{bertsekas1999nonlinear}).
\qed

Note that the hull relaxation~\eqref{model:hr_disaggregated_convex_gdp_set_based} is a convex optimization problem, which means that Lagrange multipliers are typically only guaranteed to exist at a primal optimum under some form of constraint qualification. For example, Ceria and Soares~\cite{ceria1999convex} state that if the weak Slater condition holds for every disjunct in a convex disjunctive program in DNF, then the existence of multipliers at a primal optimum  is guaranteed. 
For our purposes, we are more concerned with how to use multipliers if they exist than with the conditions guaranteeing their existence.
}

Let $\mc{P} = \{ \mc{J}_1 , \ldots, \mc{J}_{P} \}$ be a partition of $\mc{K}$,
i.e., the sets $\mc{J}_p$ are non-empty disjoint subsets of $\mc{K}$ such that $\cup_{p=1}^{P} \mc{J}_p = \mc{K}$.
For any partition $\mc{P}$ of $\mc{K}$ and any set of vectors $\vmu_1,\ldots,\vmu_{P}$, we define the \textit{partition relaxation} as 
\begin{equation} \label{model:partition_relaxation}
L_{\mc{P}}(\vmu_1,\ldots,\vmu_{P}) =
\sum_{p=1}^{P} \min \left\{ \vmu_p^\top \v{v} : \v{v} \in \cap_{k \in \mc{J}_p} \mc{F}_k \right\}~.
\end{equation}
Each minimization problem {\color{black} (or \textit{subproblem})} in the summation of \eqref{model:partition_relaxation} is a convex disjunctive program and can be solved as a convex MINLP{\color{black}, e.g., using a formulation akin to \eqref{model:hr_disaggregated_convex_gdp_algebraic} or with a Big-M formulation as described in \cite{grossmann2013systematic}}. The next proposition states that, given any partition $\mc{P}$ of $\mc{K}$, one can simply sum the associated multipliers and expect a bound improvement{\color{black}, or at least no deterioration,} relative to the previous relaxation given by \eqref{model:lr_disaggregated_convex_gdp_set_based}. 

\begin{proposition} \label{prop:partition_relaxation_bound} 
For any set of vectors $\vlambda_1,\ldots,\vlambda_{K}$ and any partition $\mc{P}$ of $\mc{K}$, we have
\begin{equation}
L_{\mc{P}} \left( \sum_{j \in \mc{J}_1} \vlambda_j, \ldots, \sum_{j \in \mc{J}_{P}} \vlambda_j \right) 
\geq 
LR_{\mc{K}}(\vlambda_1,\ldots,\vlambda_{K})~.
\end{equation}
\end{proposition}
\proof
Since, for any $\mc{J} \in \mc{P}$, we have $\cap_{j \in \mc{J}} \mc{F}_{j} \subseteq \mc{F}_{k}$ for all $k \in \mc{J}$,
it follows that 
$$
%L_{\mc{P}} \left( \sum_{j \in \mc{J}_1} \vlambda_j, \ldots, \sum_{j \in \mc{J}_{|\mc{P}|}} \vlambda_j \right) = 
\sum_{p=1}^P \min \left\{ \sum_{k \in \mc{J}_p} \vlambda_k^\top \v{v} : \v{v} \in \cap_{k \in \mc{J}_p} \mc{F}_{k} \right\} 
\geq
\sum_{p=1}^P \sum_{k \in \mc{J}_p} \min \left\{ \vlambda_k^\top \v{v} : \v{v} \in \mc{F}_{k} \right\}~. 
%= LR_{\mc{K}}(\vlambda_1,\ldots,\vlambda_{|\mc{K}|})~.
$$
\qed

{\color{black}Assuming optimal multipliers to the hull relaxation~\eqref{model:hr_disaggregated_convex_gdp_set_based} exist,
Propositions~\ref{prop:weak_and_strong_duality} and \ref{prop:partition_relaxation_bound} lead} immediately to a computable bound on a simpler operation: the bound improvement from a basic step between two disjunctions $k$ and $l$.

\begin{corollary}
Let $\vlambda_1^*,\ldots,\vlambda_{K}^*$ be optimal multipliers to the hull relaxation \eqref{model:hr_disaggregated_convex_gdp_set_based} with $K$ disjunctions.
Without loss of generality, assume $k = K-1$ and $l=K$.
Then, the bound improvement $\Delta(k,l)$ from a basic step between two disjunctions $k$ and $l$ satisfies:  
\begin{equation}
\Delta(k,l) \geq
L_{ \{ \{1\},\ldots,\{K-2\},\{K-1,K\} \} } \left( \vlambda_1^*,\ldots,\vlambda_{K-2}^*,\vlambda_{K-1}^*+\vlambda_{K}^* \right) 
-
LR_{\mc{K}}(\vlambda_1^*,\ldots,\vlambda_{K}^*)
\geq 0~.
\end{equation}
\end{corollary}

Given multipliers $\vlambda_1,\ldots,\vlambda_{K}$ and a subset $\mc{J}$ of disjunctions,
we refer to the solution of the convex disjunctive program 
\begin{equation} \label{model:pseudo_basic_step}
\min \left\{ \sum_{k \in \mc{J}} \vlambda_k^\top \v{v} : \v{v} \in \cap_{k \in \mc{J}} \mc{F}_{k} \right\}
\end{equation}
as a \textit{pseudo basic step} with respect to $\vlambda_1,\ldots,\vlambda_{K}$ and $\mc{J}$.
{\color{black} Alternatively, one can see that a pseudo basic step is the solution of one subproblem in the partition relaxation~\eqref{model:partition_relaxation}.
The qualifier ``pseudo'' is meant to distinguish it from an actual basic step.
To be precise, recall that for an actual basic step involving a subset $\mc{J}$ of disjuncts, one has to perform the following operations to obtain a tighter relaxation: (a) intersect all disjunctions in $\mc{J}$; (b) bring their intersection into DNF, possibly creating a new disjunction with many disjuncts; and (c) solve the hull relaxation of the resulting convex disjunctive program.
As discussed in Section~\ref{sec:basic_steps}, these steps are computationally expensive because the resulting hull reformulation, typically resembling formulation~\eqref{model:hr_disaggregated_convex_gdp_algebraic}, may have many more decision variables after steps (a) and (b) are carried out.
A second deficiency of an actual basic step is that, in step (c), all discrete features of the original disjunctive program are relaxed as one optimizes over the convex closure of all disjunctions.  In contrast, a pseudo basic step retains some of the discrete features as it does not convexify the disjunctions in $\mc{J}$.
On the other hand, one can see that an advantage of an actual basic step is that the relaxation in step (c) has a more ``global'' view since it takes into account the interaction of all disjunctions simultaneously, not just those in $\mc{J}$. A pseudo basic step attempts to achieve this ``global'' view through multipliers (which are not required in an actual basic step), but without having to consider all disjunctions simultaneously. 
Note that we do not require a pseudo basic step to use optimal multipliers to the hull relaxation \eqref{model:hr_disaggregated_convex_gdp_set_based}.
Finally, our definition assumes that \eqref{model:pseudo_basic_step} is solved to provable optimality.
This need not be the case; one could prematurely terminate the solution of \eqref{model:pseudo_basic_step} and take the best available bound leading to a \textit{suboptimal} pseudo basic step.
}

Given a partition $\mc{P} = \{ \mc{J}_1 , \ldots, \mc{J}_{P} \}$ of $\mc{K}$, 
let $\mc{Q}=\{1,\ldots,P\}$ be the set of disjunctions, each in DNF, resulting from sequential basic steps of the elements in each $\mc{J}_q$ for $q \in \mc{Q}$.  That is, for each $q \in \mc{Q}$, $\mc{F}_q = \cap_{k \in \mc{J}_q} \mc{F}_{k}$ expressed in DNF.
Thus, we can speak of the convex disjunctive program defined by $\mc{P}$ resulting from a sequence of actual basic steps:
\begin{subequations} \label{model:partition_DNF}
\begin{alignat}{4}
\min~~& \v{c}^\top \v{x} & & \\
\st~~& \v{x} = \v{v}_q  & & \qquad \forall~q \in \mc{Q}  \quad (\vmu_q \in \Re^n) \\
& \v{v}_q \in \mc{F}_q  & & \qquad \forall~q \in \mc{Q}~. 
\end{alignat}
\end{subequations}
The following corollary places guaranteed upper and lower bounds on one or more pseudo basic steps
and, thus, summarizes where pseudo basic steps reside in the hierarchy of relaxations proposed by Balas \cite{balas:hierarchy:1985}. 
{\color{black}
\begin{corollary} \label{prop:upper_and_lower_bounds_on_improvement}
Let $\vlambda_1,\ldots,\vlambda_{K}$ be a set of Lagrange multiplier vectors such that $\sum_{k \in \mc{K}} \vlambda_k = \v{c}$.  Let $\mc{P} = \{ \mc{J}_1 , \ldots, \mc{J}_{P} \}$ be a partition of $\mc{K}$ and $\mc{Q}=\{1,\ldots,P\}$ be the set of disjunctions, each in DNF, resulting from sequential basic steps of the elements in each $\mc{J}_q$ for $q \in \mc{Q}$.  Then, we have
\begin{equation}
z^* 
=
L_{ \{\mc{K}\} } \left( \sum_{k \in \mc{K}} \vlambda_k \right) 
\geq
z_{\mc{Q}}^{\Hrel} 
\geq
z_{\mc{Q}}^{\LD}
\geq  
L_{\mc{P}} \left( \sum_{j \in \mc{J}_1} \vlambda_j, \ldots, \sum_{j \in \mc{J}_{P}} \vlambda_j \right) 
\geq 
LR_{\mc{K}}(\vlambda_1,\ldots,\vlambda_{K})~,
\end{equation}
with $z_{\mc{Q}}^{\Hrel} = z_{\mc{Q}}^{\LD}$ if optimal multipliers 
$\vmu_1^*,\ldots,\vmu_{P}^*$ exist to the hull relaxation of \eqref{model:partition_DNF}.
\end{corollary}
\proof
The equality follows from the definition of the partition relaxation \eqref{model:partition_relaxation} 
and the assumption $\sum_{k \in \mc{K}} \vlambda_k = \v{c}$, and states the obvious: 
the optimal value $z^*$ of \eqref{model:disaggregated_convex_gdp_set_based} is equal to the trivial partition relaxation when $\mc{P} = \{\mc{K}\}$.  
The first inequality follows from the fact that any hull relaxation of the original problem, regardless of the number of basic steps applied, is still a relaxation.
%\eqref{model:partition_DNF} is a relaxation of \eqref{model:disaggregated_convex_gdp_set_based}.
The second inequality and the final assertion on conditions when $z_{\mc{Q}}^{\Hrel} = z_{\mc{Q}}^{\LD}$ follow from Proposition~\ref{prop:weak_and_strong_duality}. 
The third inequality is due to the fact that the multipliers $\left( \sum_{j \in \mc{J}_1} \vlambda_j, \ldots, \sum_{j \in \mc{J}_{P}} \vlambda_j \right)$ are suboptimal to the Lagrangian dual problem associated with  \eqref{model:partition_DNF}.
%is due to the fact that a pseudo basic step is never better than an actual basic step. 
The last inequality is from Proposition \ref{prop:partition_relaxation_bound}.
 \qed
 
%\begin{remark}
The fact that 
$
z_{\mc{Q}}^{\Hrel} 
\geq  
L_{\mc{P}} \left( \sum_{j \in \mc{J}_1} \vlambda_j, \ldots, \sum_{j \in \mc{J}_{P}} \vlambda_j \right) 
$
in Corollary~\ref{prop:upper_and_lower_bounds_on_improvement} implies that a pseudo basic step can never improve the bound more than an actual basic step, i.e., taking an actual basic step and solving the corresponding hull relaxation.
Meanwhile, since a pseudo basic step is meant to approximate an actual basic step, it is tempting to think that pseudo basic steps may always be inferior to actual basic steps in terms of bound improvement.
The fact that $
z^* 
=
L_{ \{\mc{K}\} } \left( \sum_{k \in \mc{K}} \vlambda_k \right)
$
in Corollary~\ref{prop:upper_and_lower_bounds_on_improvement},
although obvious given the definition of a partition relaxation,
reveals that in the extreme case a pseudo basic step recovers the optimal objective function value just as repeated applications of actual basic steps would do.
%\end{remark}

}

\subsection{Illustrative example} \label{sec:illustrative_example}
\begin{figure}[htbp]
\begin{center}
\includegraphics[width=1.0\textwidth]{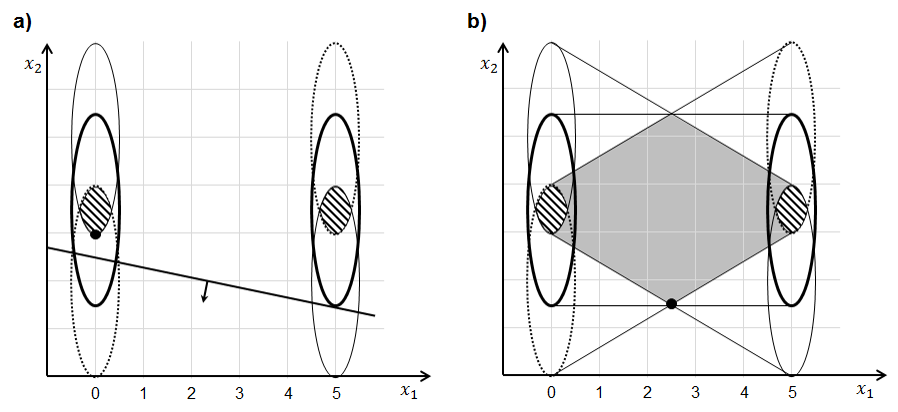}
%\vspace{-0.25in}
\caption{Illustration of example \eqref{model:simplest_example}. a) Shows the feasible region and optimal solution. b) Shows the continuous relaxation of the hull reformulation, projected onto the original space.}
\label{fig:simple_example_1}
\end{center}
\end{figure}

\begin{figure}[htbp]
\begin{center}
\includegraphics[width=1.0\textwidth]{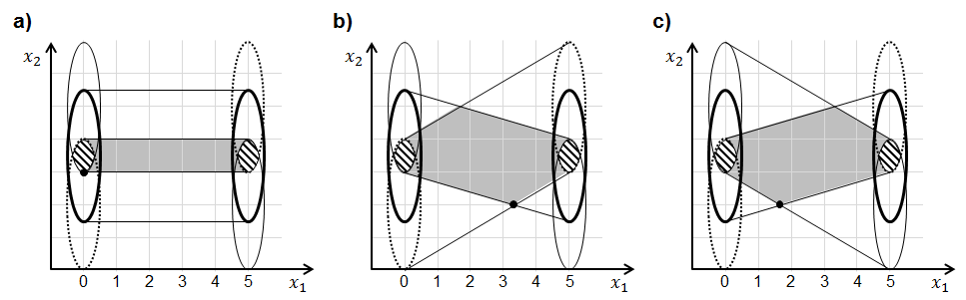}
%\vspace{-0.25in}
\caption{Illustration of the hull relaxation after basic steps for (\ref{model:simplest_example}) intersecting a) disjunctions 1 and 2; b) disjunctions 1 and 3; and c) disjunctions 2 and 3}
\label{fig:simple_example_2}
\end{center}
\end{figure}

We now present a simple convex disjunctive program in two variables to illustrate the difference in the relaxation improvement between actual and pseudo basic steps.
In light of Corollary \ref{prop:upper_and_lower_bounds_on_improvement}, which provides bounds on the improvement of a pseudo basic step, the purpose of this example is to (1) demonstrate how accurate these bounds may be, and (2) explain what goes wrong when a pseudo basic step does not achieve the same improvement as an actual basic step. Consider the following example, depicted graphically in Figure \ref{fig:simple_example_1}a:
\begin{subequations} \label{model:simplest_example}
\begin{alignat}{4}
\min~~& 0.2x_1 + x_2 & & \\
\st~~& {\color{black}\mc{F}_1 =}\left[ (x_1)^2 + \tfrac{1}{4}(x_2-5)^2 \leq 1 \right] \vee \left[ (x_1-5)^2 + \tfrac{1}{4}(x_2-2)^2 \leq 1 \right] & \\
& {\color{black}\mc{F}_2 =}\left[ (x_1)^2 + \tfrac{1}{4}(x_2-2)^2 \leq 1 \right] \vee \left[ (x_1-5)^2 + \tfrac{1}{4}(x_2-5)^2 \leq 1 \right] & \\
& {\color{black}\mc{F}_3 =}\left[ (x_1)^2 + \tfrac{1}{4}(x_2-3.5)^2 \leq 1 \right] \vee \left[ (x_1-5)^2 + \tfrac{1}{4}(x_2-3.5)^2 \leq 1 \right] & .
\end{alignat}
\end{subequations}
The first disjunction (\ref{model:simplest_example}b) is represented in Figure \ref{fig:simple_example_1}a by the ellipses with a continuous thin line, the second (\ref{model:simplest_example}c) with a dashed line, and the third (\ref{model:simplest_example}d) with a thick continuous line. The feasible region is shown in the striped area. The direction of the objective function is also shown in Figure \ref{fig:simple_example_1}a with an arrow, and the optimal solution is shown with a dot.
{\color{black}
Disaggregating the problem as in \eqref{model:disaggregated_convex_gdp_set_based} and using a compact notation, we represent problem \eqref{model:simplest_example} as follows:
\begin{subequations} \label{model:simplest_example_compact}
\begin{alignat}{4}
\min~~& 0.2x_1 + x_2 & & \\
\st~~& x_j = v_{kj} & & \qquad \forall~k \in \{1,2,3\}, j \in \{1,2\} ~~(\lambda_{kj} \in \Re) \\
& \v{v}_k \in \mc{F}_k & & \qquad \forall~k \in \{1,2,3\},
\end{alignat}
\end{subequations}
where $j$ indexes the dimensions of $\v{x}$.
The partial Lagrangian relaxation of \eqref{model:simplest_example_compact} is as follows:
\begin{subequations} \label{model:simplest_example_lagrel}
\begin{alignat}{4}
\min~~& \left(0.2 - \sum_{k \in \{1,2,3\}} \lambda_{k1} \right ) x_1 + \left(1 - \sum_{k \in \{1,2,3\}} \lambda_{k2} \right ) x_2 + \sum_{k \in \{1,2,3\}} \left( \lambda_{k1} v_{k1} + \lambda_{k2} v_{k2} \right) & & \\
\st~~& \v{v}_k \in \mc{F}_k \qquad \qquad \qquad \forall~k \in \{1,2,3\}. & &
\end{alignat}
\end{subequations}
In order to obtain the optimal multipliers that maximize \eqref{model:simplest_example_lagrel}, the hull relaxation of \eqref{model:simplest_example_compact} is solved. Note that Proposition~\ref{prop:optimal_multipliers_condition} asserts that any multipliers that maximize \eqref{model:simplest_example_lagrel} will satisfy $0.2 = \sum_{k \in \{1,2,3\}} \lambda_{k1}$ and
$1 = \sum_{k \in \{1,2,3\}} \lambda_{k2}$, while Proposition~\ref{prop:weak_and_strong_duality} implies that \eqref{model:simplest_example_lagrel} is a Lagrangian relaxation of both \eqref{model:simplest_example_compact} and its hull relaxation \eqref{model:simplest_example_hr}. The hull relaxation of \eqref{model:simplest_example} can be expressed as follows:
\begin{subequations} \label{model:simplest_example_hr}
\begin{alignat}{4}
\min~~& 0.2x_1 + x_2 & & \\
\st~~& x_j = v_{kj} & & \qquad \forall~k \in \{1,2,3\}, j \in \{1,2\} ~~(\lambda_{kj} \in \Re) \\
& \v{v}_k \in \clconv \left( \mc{F}_k \right) & & \qquad \forall~k \in \{1,2,3\},
\end{alignat}
\end{subequations}
where the convex hull of each disjunction can be expressed in algebraic form using the perspective function and the rotated second-order cone constraint described in the introduction. For example, $\clconv \left( \mc{F}_1 \right)$ of (\ref{model:simplest_example_hr}) can be represented as follows ($k=1$ below):
\begin{subequations} \label{model:simplest_example_hr_1}
\begin{alignat}{4}
&v_{k1} = \nu_{k11} + \nu_{k21} & \\
&v_{k2} = \nu_{k12} + \nu_{k22} & \\
&t_{k11} + \tfrac{1}{4}t_{k12} - \tfrac{5}{2} \nu_{k12} + \tfrac{21}{4} y_{k1} \le 0 & \\
&t_{k21} - 10 \nu_{k21} + \tfrac{1}{4}t_{k22} - \nu_{k22} + 25 y_{k2} \le 0& \\
&(\nu_{kij})^2 \le t_{kij} y_{ki} & \forall j \in \{1,2\} , i \in \{1,2\} \\
& L_{kj}y_{ki} \le \nu_{kij} \le U_{kj}y_{ki} & \forall j \in \{1,2\} , i \in \{1,2\} \\
& y_{k1} + y_{k2} = 1 & \\
& t_{kij} \ge 0 & \forall j \in \{1,2\} , i \in \{1,2\} \\
& y_{k1}, y_{k2} \ge 0~, & 
\end{alignat}
\end{subequations}
where $L_{kj}$ and $U_{kj}$ are easily computed lower and upper bounds on $\nu_{kij}$. 
A similar representation can be applied to disjunctions 2 and 3. Note that the final algebraic representation of \eqref{model:simplest_example_hr} can be more compact after algebraic substitutions. 
\begin{table}
\centering
\begin{tabular}{l c c c}
\toprule
Constraint & $\vlambda_1^*$ & $\vlambda_2^*$ & $\vlambda_3^*$ \\
\midrule
$x_1 - v_{k1}=0$ & 0.200 & 0 & 0 \\
$x_2 - v_{k2}=0$ & 0.334 & 0 & 0.666 \\
\bottomrule
\end{tabular}
\caption{Lagrange multipliers of the hull relaxation of \eqref{model:simplest_example}.}
\label{tbl:lambdas_illustrative_example}
\end{table}
}The hull relaxation of (\ref{model:simplest_example}), projected onto the original space, is presented in Figure \ref{fig:simple_example_1}b. The optimal solution is shown with a dot in Figure \ref{fig:simple_example_1}b. 
%For the hull relaxation, the optimal Lagrange multipliers are $\vlambda_1^* = [0.2,0.334]^{\top}$, $\vlambda_2^*=[0,0]^{\top}$, and $\vlambda_3^* = [0,0.666]^{\top}$ (which correspond to constraint (\ref{model:simplest_example_hr}b)). 
{\color{black}Table~\ref{tbl:lambdas_illustrative_example} lists the values of the optimal multipliers of \eqref{model:simplest_example_hr} corresponding to constraints (\ref{model:simplest_example_hr}b).} 
It is easy to see from {\color{black} these values (specifically, $\vlambda_2^* = \v{0}$)} that constraints (\ref{model:simplest_example}c), corresponding to disjunction 2, are not active. The optimal value of (\ref{model:simplest_example}) is 2.99, and the optimal value of its hull relaxation is 1.97.
{\color{black} 
The hull relaxation of the problem after applying a basic step between disjunctions 1 and 2 is as follows:
\begin{subequations} \label{model:simplest_example_hr_bs}
\begin{alignat}{4}
\min~~& 0.2x_1 + x_2 & & \\
\st~~& x_j = v_{kj} & & \qquad \forall~k \in \{1,2\}, j \in \{1,2\} ~~(\lambda_{kj} \in \Re) \\
& \v{v}_1 \in \clconv \left( \mc{F}_1 \cap \mc{F}_2 \right) & & \\
& \v{v}_2 \in \clconv \left( \mc{F}_3 \right). & & 
\end{alignat}
\end{subequations}
The problem resulting from a pseudo basic step between disjunctions 1 and 2 is as follows:
\begin{subequations} \label{model:simplest_example_lagrel_bs}
\begin{alignat}{4}
\min~~& (0.2 + 0) v_{11} + (0.334 + 0 ) v_{12} + 0 v_{21} + 0.666 v_{22} & & \\
\st~~& \v{v}_1 \in \left( \mc{F}_1 \cap \mc{F}_2 \right) & & \\
& \v{v}_2 \in \mc{F}_3~. & & 
\end{alignat}
\end{subequations}
}Table \ref{tbl:simple_example_bs} presents the objective function values of the hull relaxation after applying basic steps and pseudo basic steps. The table shows that any actual basic step and any pseudo basic step improve the hull relaxation (relative to that of the original problem, which has a value of 1.97). However, an actual basic step provides a larger improvement when applied to disjunction 1 and 2, and when applied to disjunction 2 and 3. In the first case, a basic step between disjunctions 1 and 2 yields an objective value of the relaxation that is the same as the optimal solution (which is 2.99), while the pseudo basic step gives a worse relaxation of 2.31. The illustration of the hull relaxation after basic steps, projected onto the original space, is presented in Figure \ref{fig:simple_example_2}.
\begin{table}
\centering
\begin{tabular}{c c c}
\toprule
Selected disjunctions & Basic step &Pseudo basic step \\
\midrule
1, 2 & 2.99 & 2.31 \\
1, 3 & 2.99 & 2.99 \\
2, 3 & 2.45 & 2.31 \\
\bottomrule
\end{tabular}
\caption{Objective function value of the hull relaxation of (\ref{model:simplest_example}) after intersecting disjunctions.}
\label{tbl:simple_example_bs}
\end{table}
From the initial multipliers of the hull relaxation and Figure \ref{fig:simple_example_2}a, it is easy to see why a pseudo basic step between disjunctions 1 and 2 is inferior to an actual basic step. Figure \ref{fig:simple_example_2}a shows the basic step and it is clear that disjunction~3 is inactive. The optimal Lagrange multipliers corresponding to disjunction~3 are $\vlambda_3^* = [0,0]^{\top}$, while the optimal multipliers to the new disjunction are $\vlambda_{12}^*=\v{c}=[0.2,1]^{\top}$. Meanwhile, if a pseudo basic step is applied instead of an actual basic step, the obtained (suboptimal) multipliers are $\vlambda_3 = [0,0.666]^{\top}$ and $\vlambda_{12}=[0.2,0.334]^{\top}$. {\color{black} Here, both the optimal multipliers (obtained by solving the hull relaxation of the disjunctive program after applying a basic step) and suboptimal multipliers (obtained by adding the optimal multipliers $\vlambda_k^*$ associated with the disjunctions of the original problem) satisfy $\sum_{k \in \mc{K}} \vlambda_k = \v{c}$. 
%However, the bound associated with 
However, the suboptimal multipliers provide a weaker bound than the optimal ones as shown in Table \ref{tbl:simple_example_bs} (the bound obtained using optimal multipliers is 2.99, but the one with suboptimal multipliers is 2.31).} A similar observation of the basic step between disjunctions 1 and 3 reveals that the pseudo basic step yields the same optimal multipliers as the basic step ($\vlambda_2 = \vlambda_2^* = [0,0]^{\top}$ and $\vlambda_{13} = \vlambda_{13}^* = [0.2,1]^{\top}$). Finally, the multipliers of the pseudo basic step between disjunctions 2 and 3 are $\vlambda_1 = [0.2,0.334]^{\top}$ and $\vlambda_{23} = [0,0.666]^{\top}$. In contrast, the optimal multipliers when solving the hull relaxation corresponding to the application of an actual basic step are $\vlambda_1^* = [0.3,0.501]^{\top}$ and $\vlambda_{23}^* = [-0.1,0.499]^{\top}$.

\subsection{Discussion} \label{sec:discussion}

A potential drawback of pseudo basic steps is their reliance on Lagrange multipliers.  
{\color{black} Although Proposition~\ref{prop:weak_and_strong_duality} suggests an obvious way to obtain an initial set of multipliers, i.e., by solving the hull relaxation~\eqref{model:hr_disaggregated_convex_gdp_set_based},}
for large convex disjunctive programs, it is typically more convenient computationally to remain in a lower-dimensional space using a Big-M formulation (see, e.g., \cite{grossmann2013systematic}).
As a consequence, one would see little value constructing the hull reformulation and solving its relaxation since it will not be used thereafter.
That said, it is important to emphasize that multipliers can be obtained in three ways:
1) exactly, by solving the hull relaxation to provable optimality as a convex optimization problem in polynomial time; 
2) approximately, subject to the condition {\color{black}$\sum_{k \in \mc{K}} \vlambda_k = \v{c}$}; or, in some cases,
3) analytically (see Section~\ref{sec:kmeans} below).
In fact, there is no guarantee that optimal multipliers to the hull relaxation \eqref{model:hr_disaggregated_convex_gdp_set_based} will yield greater improvements than using a set of suboptimal multipliers.
One could even use several sets of multipliers and apply pseudo basic steps using these different sets to improve the bound.

Algorithmically, the problem of choosing the disjunctions for a basic step (actual or pseudo) that will yield the most improvement could be performed in manner analogous to strong branching in mixed-integer programming.
In ``full'' strong branching, all integer variables that take a fractional value at the current node being evaluated are ``branched on'' and the one yielding the largest improvement is chosen for branching \cite{achterberg:branching:2005}.
To apply pseudo basic steps in convex disjunctive programming, one could evaluate all pairs of disjunctions, using the existing theory presented above to winnow out unproductive ones, and choose the one that yields the largest bound improvement. 
Just as full strong branching is powerful, but computationally expensive, we expect more heuristics to be developed to address this issue.

\section{Computational experiments}

In this section, we investigate two convex disjunctive programs to illustrate the improvements from pseudo basic steps and the partition relaxation. 
Although our results apply for many types of convex disjunctive programs, our computational experiments focus exclusively on MIQCP formulations to showcase the improvements possible over existing state-of-the-art solvers.
While the two problems are different in several respects, the most important distinction for our purposes is how Lagrange multipliers are obtained since our approach relies on their availability.
In Section \ref{sec:socp_dps} multipliers are computed by solving the hull relaxation, whereas in Section \ref{sec:kmeans}, we show that multipliers can be determined analytically without resorting to the hull reformulation.
All experiments were carried out on a Dell Precision T7600 running 
two 64-bit 3.10 GHz Intel Xeon E5-2687W-0 processors with 64 GB of RAM.

% Machine specs for clnl3lh342djpap
%Manufacturer: Dell
%Model: Precision T7600
%Processor: Intel Xeon CPU E5-2687W 0 at 3.10 GHz (2 processors)
%RAM: 64 GB
%64-bit Operating System 

%All experiments were carried out on a Linux machine with kernel 2.6.18 running a 64-bit x86 processor equipped with two 2.27 GHz Intel Xeon E5520 chips and 32GB of RAM.

\subsection{Conic quadratic disjunctive programs} \label{sec:socp_dps}
In this section, we numerically investigate the quality of pseudo basic steps relative to an actual basic step.
We consider a set of randomly generated conic quadratic disjunctive programs having the following form:
\begin{subequations} \label{model:socp_gdp_set_based}
\begin{alignat}{4}
\min~~& \v{c}^\top \v{x} & & \\
\st~~& \v{x} \in \bigcap_{k \in \mc{K}} \bigcup_{i \in \mc{D}_k} \mc{C}_{ki} & &
\end{alignat}
\end{subequations}
where $\mc{C}_{ki} = \left\{ \v{x} \in \Re^n : \v{x}^\top \v{Q}_{ki}\v{x} + \v{b}_{ki}^\top \v{x} + \gamma_{ki} \leq 0 \right\}$.
This section presents {\color{black}the results of two experiments with different sets of} random instances of conic quadratic disjunctive programs, each set with 10 instances. The main objective of the first experiment is to show the performance of pseudo basic steps vs. actual basic steps. This set contains small instances, since solving the hull relaxation after a few basic steps becomes intractable in larger problems. All of the instances of set 1 can be easily solved by commercial solvers. The second experiment contains larger instances that cannot be solved by CPLEX in 10 minutes. The objective in this experiment is to compare the lower bound obtained using random partitions of the problem against the lower bound of CPLEX after reaching the time limit of ten minutes.
For all instances, the variables have a lower bound of -100 and upper bound of 100. The cost coefficient of every variable is a random number between -1000 and 1000. For each disjunct $ki$ and each variable $j$, the coefficients $[q_{kijj}, b_{kij}, \gamma_{ki}]$ are random numbers between $[0.1,-1,0]$ and $[1,1,2]$.
The 10 instances of the first experiment have 8 disjunctions, with 5 disjuncts each, and $\v{x} \in \mathbb{R}^{10}$. The algorithm to randomly generate basic steps enforces that at least 4 basic steps are applied. 
\begin{table}
\centering
\begin{tabular}{cccrrcccc}
\toprule
\multicolumn{3}{c}{} & \multicolumn{2}{c}{Basic step} & \multicolumn{4}{c}{Pseudo basic step} \\
						\cmidrule(lr){4-5} 					\cmidrule(lr){6-9}  
Instance & Opt. val. & HR & relaxation & time (s) & relaxation & time (s) & parallel (s) & \% improvement \\
\midrule
1&-3,516&-4,736&-3,903&993.5&-4,295&6.5&4.9&52.9 \\
2&-3,654&-4,939&-4,054&29.7&-4,362&3.7&1.9&65.2 \\
3&-4,351&-4,868&-4,384&435.7&-4,418&5.0&2.6&93.1 \\
4&-3,470&-4,820&-3,992&1472.4&-4,321&3.9&2.3&60.2 \\
5&-5,835&-7,459&-6,453&534.2&-6,765&4.6&2.5&69.0 \\
6&-4,527&-6,542&-5,431&15.1&-5,739&6.1&2.6&72.3 \\
7&-3,377&-5,094&-3,441&154.0&-4,472&6.3&3.2&37.6 \\
8&-3,251&-4,268&-3,538&200.5&-3,628&8.2&5.7&87.8 \\
9&-3,622&-4,120&-3,976&75.6&-4,060&4.1&1.4&41.5 \\
10&-3,764&-4,872&-3,931&113.1&-4,473&5.3&2.2&42.4 \\
\bottomrule
\end{tabular}
\caption{Comparison of actual basic steps with pseudo basic steps.}
\label{tbl:random_example_bs}
\end{table}

Table~\ref{tbl:random_example_bs} presents the performance of actual basic steps vs. pseudo basic steps for these randomly generated instances.
It shows that the improvement of a pseudo basic step for these instances is between 40\% to 90\%. {\color{black} This result shows how much of the bound improvement of the basic step was captured by the pseudo basic step. The value is calculated by dividing the absolute difference between the relaxation after the pseudo basic step and the hull relaxation (column HR in Table~\ref{tbl:random_example_bs}) by the absolute difference between the relaxation after the actual basic step and the hull relaxation. The optimal multipliers for the Lagrangian dual of the problem, before basic steps, were computed by solving the hull relaxation. This is included in the solution time for the pseudo basic step results, both in serial and parallel.} 
As expected, Table~\ref{tbl:random_example_bs} shows that the relaxation after basic steps is better than the bound associated with the corresponding pseudo basic steps. However, it is important to note the required time to solve the NLP after applying multiple basic steps. For example, instance 4 requires about 25 minutes to solve the NLP and to provide a lower bound to the problem. In contrast, all of the pseudo basic steps can be computed in less than 10 seconds, and even faster when run in parallel.

The second set involves instances with 15 disjunctions, with 10 disjuncts each, and $\v{x} \in \mathbb{R}^{30}$. The algorithm to randomly generate basic steps enforces that at least 11 basic steps are applied. In this case, the randomization algorithm also enforces that at most 5 disjunctions are intersected (i.e. the resulting disjunctive program will have either 3 or 4 disjunctions. If it has 3 disjunctions, each of these will be the result of intersecting 5 of the original disjunctions. If it has 4 disjunctions, any resulting disjunction is the intersection of at most 5 of the original disjunctions). 

Table~\ref{tbl:random_example_bs_2} presents the lower bound provided by the algorithm (named TresPapas), as well as the lower bound for the hull formulation and Big-M after 10 minutes using CPLEX 12.6.3.
\begin{table}
\centering
\begin{tabular}{cccccrr}
\toprule
\multicolumn{2}{c}{} & \multicolumn{2}{c}{CPLEX (600 s)} & \multicolumn{3}{c}{TresPapas} \\
						\cmidrule(lr){3-4} 					\cmidrule(lr){5-7} 
Instance & Best known & BM LB & HR LB & LB & time (s) & parallel (s)\\
\midrule
1&-7,653&-16,933&-20,590&-14,371&234.1&110.3\\
2&-6,595&-15,498&-16,638&-10,195&265.8&257.5\\
3&-6,691&-17,731&-17,858&-12,742&136.9&95.9\\
4&-6,447&-21,127&-16,897&-12,545&111.1&69.4\\
5&-6,482&-17,071&-19,300&-12,298&328.7&176.0\\
6&-7,133&-20,026&-20,070&-12,660&210.6&102.2\\
7&-6,391&-18,289&-18,191&-12,990&214.1&73.8\\
8&-5,065&-15,215&-14,140&-9,041&160.1&109.9\\
9&-5,971&-18,202&-16,570&-10,190&150.9&67.4\\
10&-6,795&-27,926&NA&-12,246&294.9&170.1\\
\bottomrule
\end{tabular}
\caption{Comparison of lower bound for TresPapas algorithm against the lower bound of the Big-M (BM) and Hull Reformulation (HR) models after 10 minutes using CPLEX.}
\label{tbl:random_example_bs_2}
\end{table}
{\color{black} In Table \ref{tbl:random_example_bs_2}, the optimal multipliers for the Lagrangian dual of the original problems were computed by solving the hull relaxation. This is included in the solution time for the pseudo basic step results, both in serial and parallel.} Table \ref{tbl:random_example_bs_2} shows that the lower bound using TresPapas is much better than the lower bound provided by CPLEX after 10 minutes. Furthermore, 9 of the 10 instances finish in less than 5 minutes when run in serial, and in less than 2 minutes when run in parallel.

\subsection{$K$-means clustering} \label{sec:kmeans}

The $K$-means clustering problem is a widely solved problem in unsupervised learning and one of the most basic forms of clustering.  It is typically solved heuristically without any regard for a provable lower bound. In this section, we show that (1) state-of-the-art MIP solvers struggle to solve these problems to provable optimality as well as to close the optimality gap when given a near optimal solution, and (2) our partition relaxation can be much more effective at closing the gap. 

The $K$-means clustering problem is defined as follows. Suppose we are given $N$ points $\v{p}_1,\ldots,\v{p}_N$ such that $\v{p}_i \in \Re^D$, an integer $K \in \mc{N} = \{1,\dots,N\}$, and a distance metric $\D: \Re^D \times \Re^D \mapsto \Re$. 
Our goal is to determine $K$ clusters such that each point $\v{p}_i$ is assigned to a cluster with a centroid $\v{c}_k \in \Re^D$ for $k \in \mc{K} = \{1,\ldots,K\}$ and the total distance of points to centroids is minimized.  
Not only does this problem have a natural disjunctive formulation (a point must be assigned to cluster 1 or cluster 2 or ... cluster $K$), it also has a remarkably simple structure.  
It seems reasonable to expect state-of-the-art MIP solvers to be able to exploit this structure if we expect them to handle more difficult problems possessing embedded clustering decisions.  Indeed, many problems in location theory involve some sort of clustering component.  For example, the problem of determining where $K$ emergency facilities should be located to best serve a city in the wake of a natural disaster involves clustering decisions.

\subsubsection{MIQCP formulations}

When cast as a convex disjunctive program, the problem of finding $K$ optimal clusters takes the form
\begin{subequations} \label{model:Kmeans_as_convex_gdp_set_based}
\begin{alignat}{4}
\min~~& \sum_{i \in \mc{N}} d_i & & \\
\st~~& (\v{c},\v{d}) \in \bigcap_{i \in \mc{N}} \bigcup_{k \in \mc{K}} \mc{C}_{ik} & & 
\end{alignat}
\end{subequations}
where $\mc{C}_{ik} = \left\{ (\v{c},\v{d}) \in \Re^{K \times D} \times \Re_+^N : d_i \geq \D(\v{p}_i , \v{c}_k) \right\}$.
Note that there exist $N$ disjunctions, one for each point $i \in \mc{N}$, and $K$ disjuncts per disjunction.
Note that the sets $\mc{C}_{ik}$ are unbounded in this form. 
Without loss of generality, we assume that all points $\v{p}_1,\ldots,\v{p}_N$ have been normalized to reside in a $D$-dimensional hypercube and, thus, we may re-express $\mc{C}_{ik}$ as a compact set $\mc{C}_{ik} = \left\{ (\v{c},\v{d}) \in [0,1]^{K \times D} \times [0,1]^N : d_i \geq \D(\v{p}_i , \v{c}_k) \right\}$.
Henceforth, we only consider the squared $L_2$ norm as the distance metric, i.e., $\D(\v{p}_i,\v{c}_{k}) = \sum_{j=1}^D (p_{ij}-c_{kj})^2$.

%\subsubsection{$L_2$ Big M formulation}
For each point-cluster pair, let $y_{ik}$ be a binary decision variable taking value 1 if point $i$ is assigned to cluster $k$ and 0 otherwise. 
For each point $i \in \mc{N}$, define $M_i := \max\{ \D(\v{p}_i , \v{p}_{i'})  : i' \in \mc{N} \}$.
A {\color{black}straightforward} Big-M MIQCP formulation is 
\begin{subequations} \label{model:vanilla_K_means_bigM_L2}
\begin{alignat}{4}
\min_{\v{c},\v{d},\v{y}}~~& \sum_{i \in \mc{N}} d_{i} & & \\
\st~~~& d_{i} \geq \sum_{j=1}^D (p_{ij} - c_{kj})^2 - M_{i}(1-y_{ik}) & & \qquad \forall~i \in \mc{N}, k \in \mc{K} \label{eq:kmeans_bigm_distance_constraint} \\
& \sum_{k \in \mc{K}} y_{ik} = 1 & & \qquad \forall~i \in \mc{N} \label{eq:kmeans_bigm_each_point_to_exactly_one_cluster} \\
& \v{c}_k \in \Re^D & & \qquad \forall~k \in \mc{K} \\
& d_{i} \in \Re_+ & & \qquad \forall~i \in \mc{N} \\
& y_{ik} \in \{0,1\} & & \qquad \forall~i \in \mc{N}, k \in \mc{K} 
\end{alignat}
\end{subequations}
This is an extension of a typical $K$-medoids MILP formulation, which builds clusters using the $L_1$ distance metric, that one would find in an integer programming textbook (see, e.g., \cite{nemhauser:wolsey:1988}).
Clearly, if point $i$ is assigned to cluster ($y_{ik} = 1$), then constraint \eqref{eq:kmeans_bigm_distance_constraint} will be tight and $d_i$ will equal the squared Euclidean distance between point $\v{p}_i$ and centroid $\v{c}_k$.
Otherwise, this constraint is loose.  Constraints \eqref{eq:kmeans_bigm_each_point_to_exactly_one_cluster} ensure that each point is assigned to exactly one cluster. 

%\subsubsection{$L_2$ hull relaxation formulation}
Now consider a hull reformulation, also an MIQCP formulation, of the $K$-means clustering problem.
Let $M$ be a large scalar.
%The following MIQCP representation: 
\begin{subequations} \label{model:vanilla_K_means_HR_L2}
\begin{alignat}{4}
\min_{\v{c},\v{c}',\v{d},\v{d}',\v{t},\v{y}}~& \sum_{i \in \mc{N}} d_{i} & & \\
\st~~~~~~& d_{i} = \sum_{k \in \mc{K}} d'_{ik} & & \qquad \forall~i \in \mc{N} \label{eq:hr_dualized_constraint_1} \\
& c_{kj} = c'_{ki0j} + c'_{ki1j}  & & \qquad \forall~k \in \mc{K}, i \in \mc{N}, j \in \{1,\ldots,D\} \label{eq:hr_dualized_constraint_2} \\
& \sum_{j=1}^D \left( p_{ij}^2 y_{ik} - 2p_{ij}c'_{ki1j} + t_{kij} \right) - d'_{ik} \leq 0 & & \qquad \forall~i \in \mc{N}, k \in \mc{K} \label{eq:distance_copy_LB} \\
& c'^2_{ki1j} \leq t_{kij}y_{ik} & & \qquad \forall~i \in \mc{N}, k \in \mc{K}, j \in \{1,\ldots,D\} \label{eq:t_rotated_soc_constr_unit_hypercube} \\
& -M(1-y_{ik}) \leq c'_{ki0j} \leq M(1-y_{ik}) & & \qquad \forall~k \in \mc{K}, i \in \mc{N}, j \in \{1,\ldots,D\} \label{eq:CL_NO} \\
& \sum_{k \in \mc{K}} y_{ik} = 1 & & \qquad \forall~i \in \mc{N} \\
& \v{c}_k \in \Re^D & & \qquad \forall~k \in \mc{K} \\
& \v{c}'_{kib} \in [-M,M]^D & & \qquad \forall~k \in \mc{K}, i \in \mc{N}, b \in \{0,1\}  \\
& d_{i} \in \Re & & \qquad \forall~i \in \mc{N} \\
& d'_{ik} \in \Re_+ & & \qquad \forall~i \in \mc{N}, k \in \mc{K} \\
& t_{kij} \in \Re_+ & & \qquad \forall~i \in \mc{N}, k \in \mc{K}, j \in \{1,\ldots,D\} \\
& y_{ik} \in \{0,1\} & & \qquad \forall~i \in \mc{N}, k \in \mc{K} \label{eq:hr_last_constraint}
\end{alignat}
\end{subequations}

Constraints \eqref{eq:distance_copy_LB} and \eqref{eq:t_rotated_soc_constr_unit_hypercube} ensure that when point $i$ is assigned to cluster $k$, $d'_{ik}$ takes the correct nonnegative value.
Note that, in order to maintain MIQCP representability, the auxiliary decision variable $t_{kij}$ is used in constraints \eqref{eq:distance_copy_LB}, instead of the term $\frac{c'^2_{ki1j}}{y_{ik}}$,
along with the rotated second-order cone constraints \eqref{eq:t_rotated_soc_constr_unit_hypercube}. 
Constraints \eqref{eq:CL_NO} ensure that when point $i$ is not assigned to cluster $k$ ($y_{ik} = 0$), the centroid copy $\v{c}'_{ki0} = \v{0}$; otherwise, $\v{c}'_{ki0} \in [-M,M]^D$ so that $c_{kj} = c'_{ki0j} + c'_{ki1j}$ can hold for every $k$ and $i$.

%Formulation \eqref{model:vanilla_K_means_HR_L2} does not require disaggregated variables $d'_{ii'k}$ for each disjunction.
It is possible to impose tighter variable bounds, e.g., requiring $\v{c}_k \in [0,1]^D$, $\v{c}'_{kib} \in [0,1]^D$, $d_{i} \in [0,1]$, and $d'_{ik} \in [0,1]$, as well as modifying constraints \eqref{eq:CL_NO} to read $c'_{ki0j} \leq (1-y_{ik})$ (the constraint $-M(1-y_{ik}) \leq c'_{ki0j}$ vanishes when $\v{c}'_{kib} \in [0,1]^D$).
While we use the tighter bounds in computational experiments, formulation \eqref{model:vanilla_K_means_HR_L2} makes the proof in Proposition \ref{prop:Kmeans_optimal_multipliers} somewhat simpler. 

Several remarks are in order.

\begin{remark}
Both the Big-M formulation \eqref{model:vanilla_K_means_bigM_L2} and hull reformulation \eqref{model:vanilla_K_means_HR_L2} are expressed with a so-called unary representation.
One may also express constraints \eqref{eq:kmeans_bigm_each_point_to_exactly_one_cluster} as SOS1 constraints or using a binary representation \cite{vielma2011modeling}.
%This unary formulation requires $N + DK + NKD$ continuous decision variables.
\end{remark}

\begin{remark}
Suppose $(\v{c}_1,\v{c}_2,\ldots,\v{c}_K) = (\v{c}_1^*,\v{c}_2^*,\ldots,\v{c}_K^*)$ in an optimal solution.
Then, any permutation $(k_1,k_2,\ldots,k_K)$ of the indices of $\mc{K}$ such that $(\v{c}_{k_1},\v{c}_{k_2},\ldots,\v{c}_{k_K}) = (\v{c}_1^*,\v{c}_2^*,\ldots,\v{c}_K^*)$ is also part of an optimal solution.
To eliminate these equivalent optima, one can add symmetry breaking constraints 
\begin{equation} \label{eq:kmeans_symmetry_breaker}
c_{k-1,1} \leq c_{k,1}, \forall~k \in \mc{K} \setminus \{1\}~,
\end{equation}
which require the centroids to be ordered in ascending order according to the first index/dimension. 
Clearly, any index/dimension could be selected. 
\end{remark}

\begin{remark}
Assuming all points are distinct, one can also add the constraints
\begin{equation} \label{eq:at_least_one_point_per_cluster}
\sum_{i \in \mc{N}} y_{ik} \geq 1 \qquad  \forall k \in \mc{K} 
\end{equation}
to ensure that each cluster contains at least one point. 
\end{remark}

\subsubsection{Partition relaxation}

Consider a partial Lagrangian relaxation of the hull reformulation \eqref{model:vanilla_K_means_HR_L2}
in which constraints \eqref{eq:hr_dualized_constraint_1} and \eqref{eq:hr_dualized_constraint_2} are dualized with multipliers $\lambda_i$ and $\vmu_{ki}$, respectively.
This leads to the Lagrangian dual problem
\begin{equation} \label{model:Kmeans_Lagrangian_dual}
\max\{ LR_{\mc{N}}(\vlambda,\vmu) : (\vlambda,\vmu) \in \Re^N \times \Re^{K \times N \times D} \}
\end{equation}
where 
\begin{subequations} \label{model:Kmeans_Lagrangian_relaxation}
\begin{alignat}{4}
LR_{\mc{N}}(\vlambda,\vmu) = \min~~& 
\sum_{i \in \mc{N}} (1-\lambda_i) d_i + \sum_{i \in \mc{N}} \lambda_i \left( \sum_{k \in \mc{K}} d'_{ik} \right) & & \notag \\
&
+ \sum_{k \in \mc{K}} \left( \v{0} -  \sum_{i \in \mc{N}} \vmu_{ik} \right)^\top \v{c}_k
+ \sum_{k \in \mc{K}} \sum_{i \in \mc{N}} \vmu_{ik}^\top \left( \sum_{b\in\{0,1\}} \v{c}'_{kib} \right)
& & \\
\st~~& \eqref{eq:distance_copy_LB}-\eqref{eq:hr_last_constraint} & & 
\end{alignat}
\end{subequations}

The next proposition shows that, for the $K$-means clustering hull reformulation \eqref{model:vanilla_K_means_HR_L2}, optimal multipliers $\vlambda^*$and $\vmu^*$ to the Lagrangian dual problem \eqref{model:Kmeans_Lagrangian_dual}, or equivalently, to the dual of the continuous relaxation of  \eqref{model:vanilla_K_means_HR_L2}, can be computed analytically. 

\begin{proposition} \label{prop:Kmeans_optimal_multipliers}
In any optimal solution to the Lagrangian dual problem \eqref{model:Kmeans_Lagrangian_dual}, we have 
$\lambda_i^* = 1$ for all $i \in \mc{N}$ and $\vmu_{ki}^* = \v{0}$ for all $k \in \mc{K}, i \in \mc{N}$. 
\end{proposition}
\proof 
Let $\sign(a) = +1$ if $a \geq 0$ and $-1$ otherwise. 
Since $d_i$ is unconstrained in the Lagrangian relaxation problem \eqref{model:Kmeans_Lagrangian_relaxation}, 
it is clear that $1 - \lambda_i^* = 0$ for all $i \in \mc{N}$,
otherwise the minimization is unbounded as $d_i^* \rightarrow -\sign(1 - \lambda_i^*)\infty$ for some $i \in \mc{N}$. 
Now suppose $\mu_{\hat{k}\hat{i}\hat{j}}^* \neq 0$ for some $(\hat{k},\hat{i},\hat{j})$. Then, the minimization \eqref{model:Kmeans_Lagrangian_relaxation} is essentially unbounded with a solution having $c'^*_{\hat{k}\hat{i}0\hat{j}} \rightarrow -\sign(\mu_{\hat{k}\hat{i}\hat{j}}^*)M$ and $y_{\hat{i}\hat{k}}^* = 0$.  In words, point $\hat{i}$ will be assigned to some other cluster/disjunct $k' \neq \hat{k}$ so that the centroid copy $\v{c}'_{\hat{k}\hat{i}0}$ associated with not selecting cluster $\hat{k}$ can go to plus or minus infinity (or, in this formulation, some large scalar $M$). 
{\color{black} Meanwhile, $y_{\hat{k}\hat{i}}^* = 0$ and the corresponding rotated cone constraint~\eqref{eq:t_rotated_soc_constr_unit_hypercube} together imply that $c'^*_{\hat{k}\hat{i}1\hat{j}}=0$.
Since we are minimizing and $d_{\hat{i}\hat{k}}' \geq 0$, it follows that $d_{\hat{i}\hat{k}}'^* = 0$ in constraint~\eqref{eq:distance_copy_LB} showing that the objective function term $\sum_{i \in \mc{N}} \lambda_i \sum_{k \in \mc{K}} d'_{ik}$ is not affected by increasing $M$. 
}
\qed

We consider two ways of constructing a partition relaxation.
%by partitioning the points and solving independent subproblems. 
%Specifically, after partitioning the set $\mc{N}$ of points (disjunctions), we solve independent (and smaller) $K$-means clustering instances. 
Algorithm \ref{algo:assign_points_to_subproblems} describes a deterministic heuristic that relies on a primal solution.
Essentially, it attempts to distribute points assigned to the same cluster in the given primal solution to different subproblems. {\color{black}Recall that a subproblem refers to one of the minimization problems in the partition relaxation~\eqref{model:partition_relaxation}.} The rationale for such point-subproblem assignments is to force points that are not in the same cluster, and thus likely to be farther in distance from one another, to be clustered together in a subproblem and consequently increase the lower bound.
Algorithm \ref{algo:randomly_assign_points_to_subproblems} describes a randomized approach to partitioning points in which each subproblem has the same number of points. We performed limited experimentation in which subproblems could have a different number of points, but did not observe much of a benefit and so restricted ourselves to the current setting.
Randomizing the partition allows us to see if there is any benefit to using a deterministic approach and to determine what might happen if the problems are solved with a naive parallel implementation.

%Deterministic approach
%For each cluster, for each point in that cluster, assign point $i$ to subproblem $s$.

\begin{algorithm}
\caption{Deterministic heuristic to assign points to subproblems}
\label{algo:assign_points_to_subproblems}
\begin{algorithmic}[1]
\REQUIRE Assume point-cluster assignments $\hat{y}_{ik}$ are given and $\mc{S} = \{1,\ldots,S\}$ is a given set of subproblems
\STATE Set $\mc{N}_s = \emptyset~\forall s \in \mc{S}$; Set $s = 1$;
\FOR{each cluster $k \in \mc{K}$ }
\FOR{each point $i : \hat{y}_{ik} = 1$}
	\STATE $\mc{N}_s = \mc{N}_s \cup \{i\}$  
	\STATE $s = s+1$; \textbf{if} $s > S$, \textbf{then} $s = 1$;
\ENDFOR
\ENDFOR
\end{algorithmic}
\end{algorithm}

%Randomized
%We fix the desired number of points per subproblem and then randomly assign points to subproblems.

\begin{algorithm}
\caption{Randomly assign points to subproblems}
\label{algo:randomly_assign_points_to_subproblems}
\begin{algorithmic}[1]
\REQUIRE Assume $\mc{S} = \{1,\ldots,S\}$ is a given set of subproblems
\REQUIRE Assume the number $N_s$ of points per subproblem is the same for all subproblems $s$
\STATE Let $U_i$ be a uniform $(0,1)$ random variable associated with point $i \in \mc{N}$
\STATE Let $\mc{N}^{\text{Sorted}}$ be the set of all points sorted in non-increasing or non-decreasing order according to $U_i$
\STATE Assign point $i \in \mc{N}^{\text{Sorted}}$ to subproblem $s = \lceil \texttt{ord}(i)/N_s \rceil$, where $\texttt{ord}(i)$ returns the ordinal position of $i \in \mc{N}^{\text{Sorted}}$ 
\end{algorithmic}
\end{algorithm}

\subsubsection{$K$-means computational experiments}

We generated nine clustering instances with 100, 500, or 1000 points in 2, 3, or 10 dimensions for a total of 27 instances. 
We then attempted to solve each instance with $K \in \{3,5,10\}$ clusters.
A warm-start solution was provided for each instance using the best solution found from calling the \texttt{kmeans} function in MATLAB with 100 replications. It is likely that the warm-start solution is optimal or near optimal. 
Our partition relaxations, called \texttt{TresPapas}, require a partition and time limit per subproblem.
We assign 20 points to each subproblem in Algorithms \ref{algo:assign_points_to_subproblems} and \ref{algo:randomly_assign_points_to_subproblems}, a number that was determined through minimal experimentation. 
%Algorithm \ref{algo:assign_points_to_subproblems} so that each subproblem has 20 points, i.e., $|\mc{N}_s| = 20$
Each subproblem is given a time limit of $\max\{60, \texttt{TOTAL\_TIME\_LIMIT}/ |\mc{N}_s| \}$ seconds.
All models and algorithms were coded in AIMMS version 4.13.1.204 and solved serially.
%All experiments were carried out on a Linux machine with kernel 2.6.18 running a 64-bit x86 processor equipped with two 2.27 GHz Intel Xeon E5520 chips and 32GB of RAM.
Each method was given a time limit of $\texttt{TOTAL\_TIME\_LIMIT} = 900$ seconds (15 minutes).

\begin{table}[h!]
\begin{center}
%\vspace{4mm}
\begin{tabular}{ll}
	\hline
	\textbf{Algorithm} &  \textbf{Description} \\
	\hline
	\texttt{BARON BM}     & BARON 15 solving Big-M Model \eqref{model:vanilla_K_means_bigM_L2}   \\
	\texttt{CPX BM}       & CPLEX 12.6.3 solving Big-M Model \eqref{model:vanilla_K_means_bigM_L2}  \\	
	\texttt{GRB BM}       & Gurobi 6.5 solving Big-M Model \eqref{model:vanilla_K_means_bigM_L2}  \\
	\texttt{CPX BM-SOS}   & CPLEX 12.6.3 solving Big-M Model \eqref{model:vanilla_K_means_bigM_L2} with constraints \eqref{eq:kmeans_bigm_each_point_to_exactly_one_cluster} implemented as SOS1 \\	
	\texttt{GRB BM-SOS}   & Gurobi 6.5 solving Big-M Model \eqref{model:vanilla_K_means_bigM_L2} with constraints \eqref{eq:kmeans_bigm_each_point_to_exactly_one_cluster} implemented as SOS1 \\
	\texttt{BARON HR}     & BARON 15 solving Hull Reformulation \eqref{model:vanilla_K_means_HR_L2}   \\
	\texttt{CPX HR}       & CPLEX 12.6.3 solving Hull Reformulation \eqref{model:vanilla_K_means_HR_L2}  \\	
	\texttt{GRB HR}       & Gurobi 6.5 solving Hull Reformulation \eqref{model:vanilla_K_means_HR_L2}  \\
	\texttt{TresPapas Det}  & Call Algorithm \ref{algo:assign_points_to_subproblems} so that $|\mc{N}_s| = 20$; Solve each subproblem with \texttt{CPX BM} \\
	\texttt{TresPapas R}  & Call Algorithm \ref{algo:randomly_assign_points_to_subproblems} ten times so that $|\mc{N}_s| = 20$; Solve each subproblem with \texttt{CPX BM}; \\
	                      & R-Mean/R-Min/R-Max = mean/min/max improvement over all 10 trials \\
	\hline
\end{tabular}
\caption{Algorithms compared.}
\label{table:algorithms_compared}
\end{center} 
\end{table}

Table \ref{table:algorithms_compared} summarizes the solvers/methods that we compared.
We compared three prominent commercial solvers - BARON 15, CPLEX 12.6.3, and Gurobi 6.5 - all of which can solve MIQCP formulations against our partition relaxation.
Each commercial solver was tasked with solving the Big-M formulation \eqref{model:vanilla_K_means_bigM_L2}, the Big-M formulation \eqref{model:vanilla_K_means_bigM_L2} with constraints \eqref{eq:kmeans_bigm_each_point_to_exactly_one_cluster} implemented as SOS1, and the hull reformulation \eqref{model:vanilla_K_means_HR_L2}.  All three of these models also include constraints \eqref{eq:kmeans_symmetry_breaker} and \eqref{eq:at_least_one_point_per_cluster}.
The same number of binary variables are present in both formulations (BM and HR); the main difference is the number of continuous variables and constraints.
Note that BARON does not handle/accept SOS constraints, hence we could not solve a BM-SOS variant for BARON.  Limited experimentation showed that the hull reformulation with SOS1 constraints was roughly the same. 
Finally, out of fairness to CPLEX, the following remark is in order:
\begin{remark} \label{remark:aimms_cplex_interface_sos_issue}
The AIMMS-CPLEX 12.6.3 interface would \textbf{not} keep the warm-start solution in memory when Constraints \eqref{eq:kmeans_bigm_each_point_to_exactly_one_cluster} were implemented with the property SOS1, hence this could unfairly handicap CPLEX as it was forced to spend time searching for a primal solution.
\end{remark}

\begin{figure}[h]
\begin{center}
\includegraphics[width=6.5in,height=4.7in]{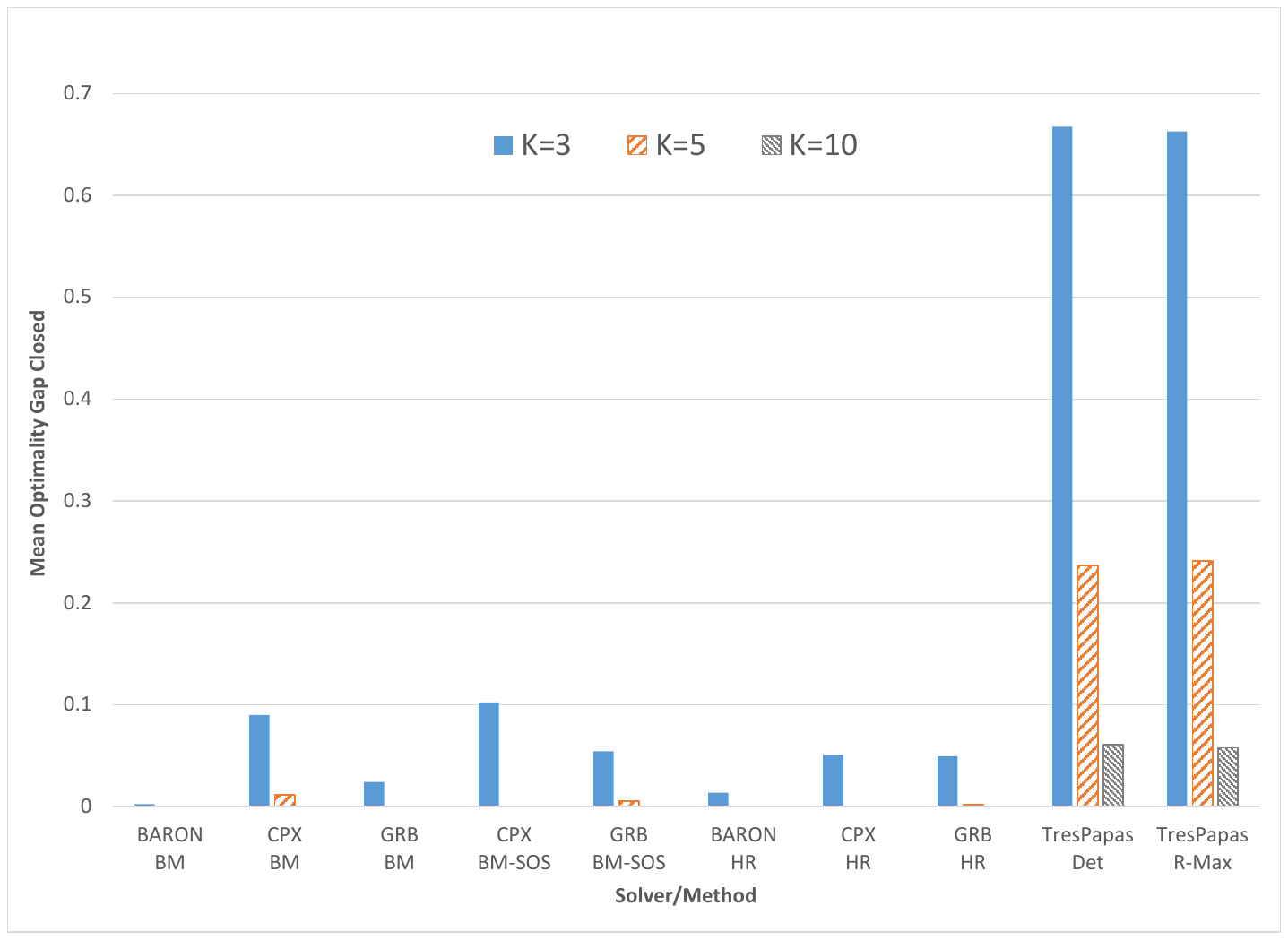}
\vspace{-0.75in}
\caption{Comparison of average optimality gap closed for various solvers and methods on $K$-means instances.}
\label{fig:kmeans_results}
\end{center}
\end{figure}

% Table generated by Excel2LaTeX from sheet 'Sheet2'
\begin{table}[h]
\footnotesize
  \centering
    \begin{tabular}{rrrrrrrrrrrrrrrrr}
    \toprule
    \multicolumn{5}{r}{}                  & \multicolumn{3}{c}{BM} & \multicolumn{2}{c}{BM-SOS} & \multicolumn{3}{c}{HR} & \multicolumn{4}{c}{TresPapas} \\ \cmidrule(lr){6-8} \cmidrule(lr){9-10} \cmidrule(lr){11-13} \cmidrule(lr){14-17}
    Instance & $N$     & $D$     & $K$     & $z^{\text{Best}}$ & \texttt{BARON} & \texttt{CPX}   & \texttt{GRB}   & \texttt{CPX}   & \texttt{GRB}   & \texttt{BARON} & \texttt{CPX}   & \texttt{GRB}   & Det   & R-Mean & R-Min & R-Max \\
\midrule
    1     & 100   & 2     & 3     & 3.14  & 0.01  & 0.35  & 0.12  & 0.40  & 0.30  & 0.06  & 0.25  & 0.24  & \textbf{0.93} & 0.79  & 0.70  & 0.85 \\
    2     & 100   & 2     & 5     & 1.91  & 0     & 0.05  & 0     & 0     & 0.05  & 0     & 0     & 0.02  & 0.61  & 0.58  & 0.50  & \textbf{0.67} \\
    3     & 100   & 2     & 10    & 0.75  & 0     & 0     & 0     & 0     & 0     & 0     & 0     & 0     & 0.04  & 0.03  & 0.02  & \textbf{0.04} \\
    4     & 100   & 3     & 3     & 4.68  & 0.01  & 0.35  & 0.09  & 0.35  & 0.19  & 0.07  & 0.19  & 0.21  & \textbf{0.87} & 0.83  & 0.81  & 0.87 \\
    5     & 100   & 3     & 5     & 3.14  & 0     & 0.05  & 0     & 0     & 0     & 0     & 0     & 0     & 0.55  & 0.49  & 0.41  & \textbf{0.55} \\
    6     & 100   & 3     & 10    & 1.79  & 0     & 0     & 0     & 0     & 0     & 0     & 0     & 0     & \textbf{0.10} & 0.08  & 0.07  & 0.10 \\
    7     & 100   & 10    & 3     & 12.68 & 0     & 0.08  & 0     & 0.1   & 0     & 0     & 0     & 0     & 0.64  & 0.61  & 0.56  & \textbf{0.66} \\
    8     & 100   & 10    & 5     & 9.70  & 0     & 0     & 0     & 0     & 0     & 0     & 0     & 0     & 0.05  & 0.05  & 0.04  & \textbf{0.05} \\
    9     & 100   & 10    & 10    & 7.86  & 0     & 0     & 0     & 0     & 0     & 0     & 0     & 0     & \textbf{0.16} & 0.14  & 0.14  & 0.16 \\
    10    & 500   & 2     & 3     & 15.02 & 0     & 0.02  & 0     & 0.01  & 0     & 0     & 0.02  & 0     & 0.75  & 0.72  & 0.67  & \textbf{0.76} \\
    11    & 500   & 2     & 5     & 8.62  & 0     & 0     & 0     & 0     & 0     & 0     & 0     & 0     & \textbf{0.40} & 0.33  & 0.29  & 0.39 \\
    12    & 500   & 2     & 10    & 4.00  & 0     & 0     & 0     & 0     & 0     & 0     & 0     & 0     & \textbf{0.03} & 0.02  & 0.02  & 0.02 \\
    13    & 500   & 3     & 3     & 17.70 & 0     & 0.01  & 0     & 0.04  & 0     & 0     & 0     & 0     & 0.81  & 0.81  & 0.79  & \textbf{0.84} \\
    14    & 500   & 3     & 5     & 13.44 & 0     & 0     & 0     & 0     & 0     & 0     & 0     & 0     & \textbf{0.17} & 0.14  & 0.13  & 0.17 \\
    15    & 500   & 3     & 10    & 8.76  & 0     & 0     & 0     & 0     & 0     & 0     & 0     & 0     & \textbf{0.04} & 0.04  & 0.03  & 0.04 \\
    16    & 500   & 10    & 3     & 67.10 & 0     & 0     & 0     & 0.01  & 0     & 0     & 0     & 0     & \textbf{0.29} & 0.26  & 0.25  & 0.28 \\
    17    & 500   & 10    & 5     & 44.40 & 0     & 0     & 0     & 0     & 0     & 0     & 0     & 0     & 0.03  & 0.03  & 0.03  & \textbf{0.03} \\
    18    & 500   & 10    & 10    & 39.51 & 0     & 0     & 0     & 0     & 0     & 0     & 0     & 0     & \textbf{0.09} & 0.08  & 0.07  & 0.08 \\
    19    & 1000  & 2     & 3     & 48.01 & 0     & 0     & 0     & 0     & 0     & 0     & 0     & 0     & \textbf{0.76} & 0.73  & 0.72  & 0.75 \\
    20    & 1000  & 2     & 5     & 25.05 & 0     & 0     & 0     & 0     & 0     & 0     & 0     & 0     & 0.19  & 0.17  & 0.15  & \textbf{0.19} \\
    21    & 1000  & 2     & 10    & 11.83 & 0     & 0     & 0     & 0     & 0     & 0     & 0     & 0     & \textbf{0.02} & 0.01  & 0.01  & 0.02 \\
    22    & 1000  & 3     & 3     & 29.21 & 0     & 0     & 0     & 0     & 0     & 0     & 0     & 0     & 0.81  & 0.79  & 0.77  & \textbf{0.81} \\
    23    & 1000  & 3     & 5     & 20.69 & 0     & 0     & 0     & 0     & 0     & 0     & 0     & 0     & \textbf{0.11} & 0.07  & 0.06  & 0.09 \\
    24    & 1000  & 3     & 10    & 14.29 & 0     & 0     & 0     & 0     & 0     & 0     & 0     & 0     & \textbf{0.02} & 0.02  & 0.01  & 0.02 \\
    25    & 1000  & 10    & 3     & 119.21 & 0     & 0     & 0     & 0     & 0     & 0     & 0     & 0     & 0.13  & 0.12  & 0.11  & \textbf{0.14} \\
    26    & 1000  & 10    & 5     & 87.56 & 0     & 0     & 0     & 0     & 0     & 0     & 0     & 0     & 0.02  & 0.01  & 0.01  & \textbf{0.02} \\
    27    & 1000  & 10    & 10    & 79.13 & 0     & 0     & 0     & 0     & 0     & 0     & 0     & 0     & 0.04  & 0.04  & 0.04  & \textbf{0.04} \\
\hline
          &       &       &       & Mean  & 0.00  & 0.03  & 0.01  & 0.03  & 0.02  & 0.00  & 0.02  & 0.02  & \textbf{0.32} & 0.30  & 0.27  & 0.32 \\
    \bottomrule
    \end{tabular}%
  \caption{Fraction of optimality gap closed in 900 s for $K$-means clustering instances.}
  \label{table:kmeans_results}%
\end{table}%

%Take-aways
Figure \ref{fig:kmeans_results} summarizes the instance-by-instance results in Table \ref{table:kmeans_results}.
{\color{black} In all instances, the MATLAB-provided warm-start solution was never improved upon and its objective function value was used as the upper bound in the optimality gap calculation.}
It is immediately clear from Figure \ref{fig:kmeans_results} that the commercial solvers struggle to close the optimality gap for even the smallest instances. Meanwhile, our partition relaxation (deterministic or randomized) is significantly better than all commercial solvers at improving the lower bound. As more clusters are included (i.e., as $K$ increases), the ability to close the gap decreases quickly. In these instances, the \texttt{TresPapas} approaches have difficulty making bound improvements, but at least make some progress in contrast with the commercial solvers which make virtually no improvements. Based on limited testing, we conjecture that the differences would be even more pronounced with a longer time limit.

When comparing the deterministic and randomized \texttt{TresPapas} approaches in Table \ref{table:kmeans_results}, it is interesting to note that the deterministic variant outperforms the randomized one on average.  Meanwhile, the randomized version appears to yield relatively stable results as can be seen by taking the difference of the R-Max and R-Min columns.  This consistency is positive news since a randomized version is remarkably easy to implement.  It would be interesting to see if the same consistency emerges in other problem classes. 

We close this section by making some higher-level remarks. 
First, we believe that the $K$-means clustering problem is a prototypical MIQCP and, while one may not expect a generic MIQCP solver to be able to solve such problems as fast as, say, knapsack problems, it seems fair to expect such solvers to be able to make modest improvements on small instances. The instances considered here are quite modest compared to what is found in the machine learning literature. 
Second, it is surprising that our approach, which is quite general as it applies to almost all convex disjunctive programs and therefore to a large number of convex MINLPs, is able to make significant improvement relative to state-of-the-art solvers.
Third, our approach is inherently parallelizable and thus we expect further gains from a parallel implementation.

\section{Conclusions and future research directions}

Basics steps are not a regularly used tool in the MIP repertoire. Nevertheless, they offer a very simple way to improve the formulation (i.e., the continuous relaxation) of a convex MINLP represented as a convex disjunctive program.  In this work, we introduced the notion of a pseudo basic step to determine guaranteed bounds on the improvement from a basic step and showed that pseudo basic steps can  
make significant bound improvements relative to state-of-the-art commercial mixed-integer programming solvers {\color{black}on several sets of conic quadratic disjunctive programming instances}. 
Given the growing excitement for solving convex MINLPs and the success of our approach, 
we hope to revive interest in basic steps and other techniques outside of the mainstream MIP toolbox.

{ \color{black}
Despite its efficacy to generate a tighter formulation, a basic step (as introduced by Balas \cite{balas:1979}) requires one to intersect two disjunctions, thus forming a new single disjunction, possibly with more disjuncts, and then to re-solve the resulting hull relaxation.
Practically speaking, a pseudo basic step retains some of the discrete features in the original disjunctive program (or ancestral disjunctive program in an iterative approach) by intersecting two or more disjunctions, and requires the solution (not necessarily to provable optimality) of a small convex disjunctive program.  Today, some readers may argue that solving a convex MINLP, however small, is prohibitively expensive.  We are more optimistic and believe that in the next two decades, our ability to solve convex MINLPs will witness roughly the same orders-of-magnitude improvement that computational MILP has experienced over the past two decades.  In fact, we are hopeful that the adage ``Just MIP it!'' \cite{fischetti2009justmipit} will soon apply to convex MINLPs as easily as it applies to MILPs today. 
}

There are several fruitful pathways for additional research.
Algorithmically, the question remains: How should basic steps be optimally chosen?  Just as it is difficult to determine the best branching rule or cut in MIP, we conjecture that this question does not have a trivial answer.  Nevertheless, there are likely some guiding principles that could help in computations. Another open question is: How can cuts generated from pseudo basic steps be included in a branch-and-cut framework?
%Is there a systematic way to choose a partition for the partition relaxation?
As far as applications are concerned, since our approach was successful when focused exclusively on MIQCP formulations, we hope to explore if there is even greater potential when applied to more general convex MINLP problems.
Theoretically, it would be interesting to better understand the correspondence between the convex disjunctive program \eqref{model:pseudo_basic_step} associated with a pseudo basic step and the various cut generating optimization problems used, for example, to generate lift-and-project cuts.
In addition, since much of the cutting-plane theory associated with linear disjunctive programming relies on two-term disjunctions, with some notable exceptions, e.g., Perregaard and Balas \cite{perregaard:balas:2001} consider lift-and-project cuts from multi-term linear disjunctions, it may be promising to investigate cut generating mechanisms for multi-term convex disjunctions given the success of our approach.

\section*{Acknowledgments}

We wish to thank Ignacio Grossmann, Nick Sawaya, Myun-Seok Cheon, and Ahmet Keha for their feedback on a preliminary manuscript.
{\color{black}We are grateful to two anonymous referees whose perceptive comments significantly improved the quality of the paper.}

\small
\bibliographystyle{abbrv}
\bibliography{gdp_refs,reference,branching_refs}

\begin{thebibliography}{10}

\bibitem{achterberg:branching:2005}
T.~Achterberg, T.~Koch, and A.~Martin.
\newblock Branching rules revisited.
\newblock {\em Operations Research Letters}, 33(1):42--54, 2005.

\bibitem{balas:1979}
E.~Balas.
\newblock Disjunctive programming.
\newblock {\em Annals of Discrete Mathematics}, 5:3--51, 1979.

\bibitem{balas:hierarchy:1985}
E.~Balas.
\newblock Disjunctive programming and a hierarchy of relaxations for discrete
  optimization problems.
\newblock {\em SIAM Journal on Algebraic Discrete Methods}, 6(3):466--486,
  1985.

\bibitem{balas1998disjunctive}
E.~Balas.
\newblock Disjunctive programming: Properties of the convex hull of feasible
  points.
\newblock {\em Discrete Applied Mathematics}, 89(1):3--44, 1998.

\bibitem{balas:2010}
E.~Balas.
\newblock {\em Disjunctive Programming}, pages 283--340.
\newblock Springer Berlin Heidelberg, Berlin, Heidelberg, 2010.

\bibitem{balas:lift_and_project:1993}
E.~Balas, S.~Ceria, and G.~Cornu{\'e}jols.
\newblock A lift-and-project cutting plane algorithm for mixed 0--1 programs.
\newblock {\em Mathematical programming}, 58(1-3):295--324, 1993.

\bibitem{belotti:2016}
P.~Belotti, P.~Bonami, M.~Fischetti, A.~Lodi, M.~Monaci, A.~Nogales-G{\'o}mez,
  and D.~Salvagnin.
\newblock On handling indicator constraints in mixed integer programming.
\newblock {\em Computational Optimization and Applications}, pages 1--22, 2016.

\bibitem{bental2001lectures}
A.~Ben-Tal and A.~Nemirovski.
\newblock {\em Lectures on modern convex optimization: analysis, algorithms,
  and engineering applications}.
\newblock SIAM, 2001.

\bibitem{bertsekas1999nonlinear}
D.~P. Bertsekas.
\newblock {\em Nonlinear programming}.
\newblock Athena Scientific, 2nd edition, 1999.

\bibitem{bonami:2015}
P.~Bonami, A.~Lodi, A.~Tramontani, and S.~Wiese.
\newblock On mathematical programming with indicator constraints.
\newblock {\em Mathematical Programming}, 151(1):191--223, 2015.

\bibitem{ceria1999convex}
S.~Ceria and J.~Soares.
\newblock Convex programming for disjunctive convex optimization.
\newblock {\em Mathematical Programming}, 86(3):595--614, 1999.

\bibitem{fischetti2009justmipit}
M.~Fischetti, A.~Lodi, and D.~Salvagnin.
\newblock Just {MIP} it!
\newblock In {\em Matheuristics}, pages 39--70. Springer, 2009.

\bibitem{grossmann2013systematic}
I.~E. Grossmann and F.~Trespalacios.
\newblock Systematic modeling of discrete-continuous optimization models
  through generalized disjunctive programming.
\newblock {\em AIChE Journal}, 59(9):3276--3295, 2013.

\bibitem{guignard2003lagrangean}
M.~Guignard.
\newblock Lagrangean relaxation.
\newblock {\em Top}, 11(2):151--200, 2003.

\bibitem{nemhauser:wolsey:1988}
G.~L. Nemhauser and L.~A. Wolsey.
\newblock {\em Integer and combinatorial optimization}.
\newblock John Wiley \& Sons, 1988.

\bibitem{perregaard:balas:2001}
M.~Perregaard and E.~Balas.
\newblock {\em Generating Cuts from Multiple-Term Disjunctions}, pages
  348--360.
\newblock Springer Berlin Heidelberg, Berlin, Heidelberg, 2001.

\bibitem{raman1994modelling}
R.~Raman and I.~E. Grossmann.
\newblock Modelling and computational techniques for logic based integer
  programming.
\newblock {\em Computers \& Chemical Engineering}, 18(7):563--578, 1994.

\bibitem{ruiz:hierarchy:2012}
J.~P. Ruiz and I.~E. Grossmann.
\newblock A hierarchy of relaxations for nonlinear convex generalized
  disjunctive programming.
\newblock {\em European Journal of Operational Research}, 218(1):38--47, 2012.

\bibitem{sawaya2006reformulations}
N.~Sawaya.
\newblock {\em Reformulations, relaxations and cutting planes for generalized
  disjunctive programming}.
\newblock PhD thesis, 2006.

\bibitem{sawaya:hierarchy:2012}
N.~Sawaya and I.~Grossmann.
\newblock A hierarchy of relaxations for linear generalized disjunctive
  programming.
\newblock {\em European Journal of Operational Research}, 216(1):70--82, 2012.

\bibitem{stubbs1999branch}
R.~A. Stubbs and S.~Mehrotra.
\newblock A branch-and-cut method for 0-1 mixed convex programming.
\newblock {\em Mathematical programming}, 86(3):515--532, 1999.

\bibitem{trespalacios2014algorithmic}
F.~Trespalacios and I.~E. Grossmann.
\newblock Algorithmic approach for improved mixed-integer reformulations of
  convex generalized disjunctive programs.
\newblock {\em INFORMS Journal on Computing}, 27(1):59--74, 2014.

\bibitem{trespalacios:lagrangean:2016}
F.~Trespalacios and I.~E. Grossmann.
\newblock Lagrangean relaxation of the hull-reformulation of linear generalized
  disjunctive programs and its use in disjunctive branch and bound.
\newblock {\em European Journal of Operational Research}, 253(2):314 -- 327,
  2016.

\bibitem{vielma2011modeling}
J.~P. Vielma and G.~L. Nemhauser.
\newblock Modeling disjunctive constraints with a logarithmic number of binary
  variables and constraints.
\newblock {\em Mathematical Programming}, 128(1-2):49--72, 2011.

\end{thebibliography}

\end{document}